\documentclass[12pt,leqno]{amsart}
\usepackage{amssymb}
\usepackage{amsmath,amssymb}
\renewcommand{\baselinestretch}{1.1}

\newtheorem{theorem}{Theorem}[section]

\newtheorem{proposition}{Proposition}[section]

\begin{document}
\theoremstyle{plain}
\newtheorem{MainThm}{Theorem}
\newtheorem{thm}{Theorem}[section]
\newtheorem{clry}[thm]{Corollary}
\newtheorem{prop}[thm]{Proposition}
\newtheorem{lem}[thm]{Lemma}
\newtheorem{deft}[thm]{Definition}
\newtheorem{hyp}{Assumption}
\newtheorem*{ThmLeU}{Theorem (J.~Lee, G.~Uhlmann)}

\theoremstyle{definition}
\newtheorem{rem}[thm]{Remark}
\newtheorem*{acknow}{Acknowledgments}
\numberwithin{equation}{section}
\newcommand{\eps}{{\varphi}repsilon}
\renewcommand{\d}{\partial}
\newcommand{\re}{\mathop{\rm Re} }
\newcommand{\im}{\mathop{\rm Im}}
\newcommand{\R}{\mathbf{R}}
\newcommand{\C}{\mathbf{C}}
\newcommand{\N}{\mathbf{N}}
\newcommand{\D}{C^{\infty}_0}
\renewcommand{\O}{\mathcal{O}}
\newcommand{\dbar}{\overline{\d}}
\newcommand{\supp}{\mathop{\rm supp}}
\newcommand{\abs}[1]{\lvert #1 \rvert}
\newcommand{\csubset}{\Subset}
\newcommand{\detg}{\lvert g \rvert}
\newcommand{\dd}{\mbox{div}\thinspace}
\newcommand{\www}{\widetilde}
\newcommand{\ggggg}{\mbox{\bf g}}
\newcommand{\ep}{\varepsilon}
\newcommand{\la}{\lambda}
\newcommand{\va}{\varphi}
\newcommand{\ppp}{\partial}
\newcommand{\ooo}{\overline}
\newcommand{\wwwx}{\widetilde{x}_0}
\newcommand{\sumkj}{\sum_{k,j=1}^n}
\newcommand{\walpha}{\widetilde{\alpha}}
\newcommand{\wbeta}{\widetilde{\beta}}
\newcommand{\weight}{e^{2s\va}}
\newcommand{\fdif}{\partial_t^{\alpha}}
\newcommand{\OOO}{\Omega}
\newcommand{\LLL}{L_{\lambda,\mu}}
\newcommand{\LLLL}{L_{\widetilde{\la},\widetilde{\mu}}}
\newcommand{\ppdif}[2]{\frac{\partial^2 #1}{{\partial #2}^2}}
\newcommand{\uu}{\mbox{\bf u}}
\newcommand{\vv}{\mbox{\bf v}}
\newcommand{\ww}{\mbox{\bf w}}
\renewcommand{\v}{\mathbf{v}}
\newcommand{\y}{\mathbf{y}}
\newcommand{\ddd}{\mbox{div}\thinspace}
\newcommand{\rrr}{\mbox{rot}\thinspace}
\newcommand{\Y}{\mathbf{Y}}
\newcommand{\w}{\mathbf{w}}
\newcommand{\z}{\mathbf{z}}
\newcommand{\G}{\mathbf{G}}
\newcommand{\f}{\mathbf{f}}
\newcommand{\F}{\mathbf{F}}
\newcommand{\dddx}{\frac{d}{dx_0}}
\newcommand{\CC}{_{0}C^{\infty}(0,T)}
\newcommand{\HH}{_{0}H^{\alpha}(0,T)}
\newcommand{\llll}{L^{\infty}(\Omega\times (0,t_1))}
\renewcommand{\baselinestretch}{1.5}
\renewcommand{\div}{\mathrm{div}\,}  
\newcommand{\grad}{\mathrm{grad}\,}  
\newcommand{\rot}{\mathrm{rot}\,}  

\title
[]
{Remark on controllability to trajectories of  a simplified fluid-structure iteration model.}


\author{
O.~Yu.~Imanuvilov }

\thanks{
Department of Mathematics, Colorado State
University, 101 Weber Building, Fort Collins, CO 80523-1874, U.S.A.
e-mail: {\tt oleg@math.colostate.edu}
Partially supported by NSF grant DMS 1312900}\,


\date{}

\maketitle

\begin{abstract}
We prove the exact  controllability  result to trajectories  of a simplified model of motion of  a rigid body in  fluid  flow. Unlike a previously know  results such  a trajectory does not need to be a stationary solution.
\end{abstract}

\section{Introduction and main results}

The paper is concern with the following controllability problem: In the bounded domain $Q=(0,T)\times \Omega, \Omega=[a,b]$,  $-\infty< a<0<b<+\infty, x=(x_0,x_1)$ we consider the system of semilinear heat equations
\begin{eqnarray}\label{eq1}
G_1(x,w_1)= \rho_1\partial_{x_0} w_1- a_1\partial_{x_1}^2w_1+b_1\partial_{x_1}w_1+c_1w_1\\ +g_1(x,w_1,\partial_{x_1}w_1)=f_1+\chi_\omega u\quad \mbox{in}\,\, Q_+=(0,T)\times (0,b),\nonumber
\end{eqnarray}
\begin{eqnarray}
\label{e2}
G_2(x,w_2)=\rho_2\partial_{x_0} w_2-a_2\partial_{x_1}^2w_2+b_2\partial_{x_1}w_2+c_2w_2\\ +g_2(x,w_2,\partial_{x_1}w_2)=f_2\quad \mbox{in}\,\, Q_-=(0,T)\times (a,0).\nonumber
\end{eqnarray} On the interface $[0,T]\times \{0\}$ functions $w_1,w_2$ are connected through the boundary conditions
\begin{equation}\label{eq3}
w_1(x_0,0)-w_2(x_0,0)=(\partial_{x_1}w_1-\partial_{x_1}w_2-M\partial_{x_0}w_1)(x_0,0)+ r(x_0) \quad \mbox{on} \,\,[0,T],
\end{equation}
where  $r$ is a given function, $M$ is a positive constant. On the  lateral boundary of cylinder $Q$ functions  $w_1,w_2$ satisfy the Dirichlet boundary conditions
\begin{equation}\label{eeq4}
w_1(x_0,b)=w_2(x_0,a)=0 \quad \mbox{on} \,\,[0,T].
\end{equation}
The initial condition is
\begin{equation}\label{eq4}
w(0,x_1)=w_0(x_1)\quad \mbox{on}\,\,\Omega.
\end{equation}
Here
$$
w=\left \{ \begin{matrix} w_1\quad \mbox{for}\quad x\in Q_+,\\ w_2\quad \mbox{for} \quad x\in Q_-.\end{matrix}\right.
$$
Function $u(x)$ is the control distributed over domain $Q_\omega=(0,T)\times \omega,  \omega=(d,b), d\in (0,b): \mbox{supp}\, u\subset Q_\omega.$
One of the physical applications of system (\ref{eq1})-(\ref{eeq4}) is the rigid body moving through the fluid flow, where  $w_j$ is velocity of the fluid flow, $M$ is the mass of the rigid body, $\int_0^{x_0}w_1(\tilde x_0,0)d\tilde x_0 +h_0$ is the position of the body (see \cite{VZ} for details of the model.)  We are looking for locally distributed  control $u$ such that at moment $T$  for the given target function $w_2$ we have:
\begin{equation}\label{eq5}
w(T,x_1)=w_2(x_1)\quad \mbox{on}\,\,\Omega.
\end{equation}
We make the following  standard assumptions:
\begin{equation}\label{eq6}
\rho_1, a_1, b_1\in C^1(\bar Q_+), \,\, \rho_2, a_2, b_2\in C^1(\bar Q_-),\quad  c_1\in  L^\infty(0,T;L^2(0,b)),\quad c_2\in L^\infty(0,T;L^2(a,0)),
\end{equation}
there exist a positive constant $\alpha$ such that
\begin{equation}\label{eq7}
\rho_1(x)\ge \alpha>0,\quad
a_1(x)\ge \alpha>0 \quad \forall x\in Q_+,\quad \rho_2(x)\ge \alpha>0,\quad  a_2(x)\ge \alpha>0 \quad \forall x\in\, Q_-,
\end{equation}
\begin{equation}\label{eq8}
g_1\in C^2(\bar Q_+\times \Bbb R^1\times \Bbb R^1),\quad g_2\in C^2(\bar Q_-\times \Bbb R^1\times \Bbb R^1).
\end{equation} there exist constants $C_1,\dots , C_2 $ independent of  $x$ and $\xi_i$, and  $p_j\ge 1, j\in\{1,2,3\}$ such that
\begin{eqnarray}\label{eq9}
\vert g_i(x,\xi_1,\xi_2)\vert\le C(1+\vert \xi_1\vert^{p_1}+\vert \xi_1\vert^{p_1}\vert \xi_2\vert), \quad \vert \partial_{\xi_1}g_i(x,\xi_1,\xi_2)\vert\le C(1+\vert \xi_1\vert^{p_1-1}+\vert \xi_1\vert^{p_2-1}\vert \xi_2\vert),\nonumber\\\quad
\vert \partial_{\xi_2}g_i(x,\xi_1,\xi_2)\vert\le C(1+\vert \xi_1\vert^{p_3})\quad \forall (x,\xi_1,\xi_2)\in Q\times \Bbb R^2\quad\mbox{and}\quad \forall i\in \{1,2\}.
\end{eqnarray}

{\bf Remark.} {\it  The nonlinear term $g_1(x,u,\partial_{x_1}u)=g_2(x,u,\partial_{x_1}u)=u\partial_{x_1} u$ satisfies (\ref{eq8}) and (\ref{eq9}).}

Since it is well known that the controllability problem (\ref{eq1})-(\ref{eq5}) can not be solved for an arbitrary  target function $w_2$ we introduce the additional condition:

{\bf Condition 1.} {\it  There exist a pair $ \mbox{\bf w}=(\mbox{\bf w}_1,\mbox{\bf w}_2)\in H^{1,2}(Q_+)\times H^{1,2}(Q_-)$ and control $\mbox{\bf u}\in L^2(Q), \mbox{supp}\, \mbox{\bf u}\subset \bar Q_\omega$ such that
\begin{equation}\label{Eeq1}\nonumber
(G_1(x,\mbox{\bf w}_1),G_2(x,\mbox{\bf w}_2))=(f_1+\mbox{\bf u},f_2),
\end{equation}
\begin{equation}\label{Ee3}\nonumber
(\mbox{\bf w}_1-\mbox{\bf w}_2)(x_0,0)=(\partial_{x_1}\mbox{\bf w}_1-\partial_{x_1}\mbox{\bf w}_2-M\partial_{x_0}\mbox{\bf w}_1)(x_0,0)- r(x_0)=0 \quad \mbox{on} \,\,\,[0,T],
\end{equation}
\begin{equation}\nonumber\mbox{\bf w}(\cdot,b)=\mbox{\bf w}(\cdot,a)=0 \quad \mbox{on} \,\,\,[0,T],
\quad
\mbox{\bf w}(T,\cdot)=w_2.
\end{equation}}
Our main  result is the following
\begin{theorem}\label{gonduras}  Let $f_1\in L^2(Q_+),f_2\in L^2(Q_-), r\in L^2(0,T).$ Suppose that assumptions (\ref{eq5})- (\ref{eq9}) and  Condition 1  for functions  $f_1,f_2,r$  holds true, $w_0\in H^1_0(\Omega)$. Then there exists a positive $\epsilon>0$ such that if
$$
\Vert w_0-\mbox{\bf w}(\cdot,0)\Vert_{H^1_0(\Omega)}\le \epsilon$$ the controllability problem (\ref{eq1})-(\ref{eq5}) has solution $(w,u)\in H^{1,2}(Q_+)\cap H^{1,2}(Q_-)\cap C^0(0,T; H^1_0(\Omega))\times L^2(Q), \mbox{supp}\, u\subset \bar Q_\omega .$
\end{theorem}
Theorem \ref{gonduras} was established for the case $\mbox{\bf w}\equiv 0$  with  control located at both ends in \cite {D-F} and at one end  in \cite{LTT}.

Another physical application of the controllability problem (\ref{eq1})-(\ref{eq5}) describes to rods connected by a point mass. (see \cite {HM1} for details of the  model.) The zero null controllability for this model was proved in \cite {HM2} for the case when coefficients  $\rho_j, a_j, b_, c_j$ are constants and recently in \cite{AB} when  coefficients  $\rho_j, a_j,b_, c_j$  are space dependent functions. The method of both papers based on the analysis of eigenvalues and eigenfunctions and therefore can  not be applied to the case  of time dependent coefficients.

The n-dimensional generalization  for linear  parabolic equations with time independent coefficients  was studied by J.L. Russeau with co-authors in  \cite{Russeau}. Exact  controllability  of similar problem for linear  1-d hyperbolic equations  in case of one point mass attached was proved  by S. Hansen and E.Zuazua in \cite{Scott} and for several mass attached  case by     S. Avdonin and J. Edwards in  \cite{AE2018}.

The proof of the Theorem \ref{gonduras} is based on the  implicit  function theorem and  null-controllability result  for the linearized system (\ref{eq1}) - (\ref{eeq4}). The null-controllability of the linearized system follows from the observability estimate. Observability estimate is   obtained by Carleman estimate with boundary. The weight function is similar to one from work \cite{Dubova}.

{\bf Notations.} Let $\Omega_+=(0,b), \Omega_-=(a,0),$ $i = \sqrt{-1}$ and $D= (D_0, D_1)$, $D_0=\frac{1}{i}
\partial_{x_0},D_1=\frac{1}{i}
\partial_{x_1}, \alpha=(\alpha_0,\alpha_1), \alpha_0\ge 0,\alpha_1\ge 0, \vert\alpha\vert=2\alpha_0+\alpha_1, \partial^\alpha=\partial^{\alpha_0}_{x_0}\partial^{\alpha_1}_{x_1}.$ For any function $\tilde \rho$ we introduce the space $L^2_{\tilde \rho}(X)=\{u \vert \Vert u\Vert_{L^2_{\tilde\rho}(X)}=\root\of{\int_X \vert\tilde\rho\vert u^2dx}\},$ by $Fu$ we denote the Fourier transform of function $u$ in variable $x_0$: $Fu=\frac{1}{\root\of{2\pi}}\int_{-\infty}^\infty e^{-i\xi_0x_0} u(x_0)dx_0.$ Let $\xi=(\xi_0,\xi_1),
\zeta=(\xi_0,\tilde s),$ $x=(x_0,x_1),$ $ \zeta^*=(\xi_0^*,\tilde s^*), M(\xi_0,\tilde s)=(\tilde s^4+\xi_0^2)^\frac 14, \Bbb M=\{( \xi_0,\widetilde s); \thinspace
M( \xi_0, \widetilde s)=1\}.$ We introduce the conic neighborhood of the point $\zeta^*$:
$$
\mathcal O(\zeta^*,\delta)=\{ (\xi_0,\tilde s)\in \Bbb R^2\setminus \{0\} \vert (\xi_0/M^2(\xi_0,\tilde s),\tilde s/M(\xi_0,\tilde s))-(\xi_0^*,\tilde s^*)\vert\le \delta\},
$$
and the Sobolev spaces
$$
H^{1,2}(Q_\pm)=\{ u\vert u, \partial_{x_0}u, \partial_{x_1} u, \partial_{x_1}^2u\in L^2(Q_\pm)\},
$$
$ H^{1,2,\tilde s}(Q_\pm)$ the space $H^{1,2}(Q_\pm)$ equipped with the norm
$$
\Vert u\Vert_{H^{1,2,\tilde s}(Q_\pm)}=\root\of{ \Vert u\Vert^2_{H^{1,2}(Q_\pm )}  +\tilde s^2\Vert u\Vert^2_{L^2(Q_\pm)} }.
$$
For any function $u$ we set $[u]=\lim_{x_1\rightarrow +0} u(x_0,x_1)-\lim_{x_1\rightarrow -0} u(x_0,x_1).$ For the symbol $M(\xi_0,\tilde s)$ we introduce the pseudodifferential operator by formula
$$
M(D_0,\tilde s)u=\frac{1}{\root\of{2\pi}}\int_{\Bbb R^1} M(\xi_0,\tilde s)e^{i\xi_0 x_0} F u d\xi_0.
$$

\section{Observability Estimate.}

 In this section we prove the observability estimate for the following system:
\begin{equation}\label{Beq1}
 -\tilde \rho_1\partial_{x_0} v_1-\tilde a_1\partial_{x_1}^2v_1+\tilde b_1\partial_{x_1}v_1+\tilde c_1v_1=\tilde f_1\quad \mbox{in}\,\, Q_+,
\end{equation}
\begin{equation}
\label{Be2}
-\tilde \rho_2\partial_{x_0} v_2- \tilde a_2\partial_{x_1}^2v_2+\tilde b_2\partial_{x_1}v_2+\tilde c_2v_2=\tilde f_2\quad \mbox{in}\,\, Q_-.
\end{equation} On the interface $[0,T]\times \{0\}$ functions $v_1,v_2$ are connected through the boundary conditions
\begin{equation}\label{BBe3}
v_1(x_0,0)-v_2(x_0,0)=(\partial_{x_1}v_1-\partial_{x_1}v_2+M\partial_{x_0}v_1)(x_0,0)- \tilde r(x_0)=0 \quad \mbox{on}\,\, \,[0,T],
\end{equation}
\begin{equation}\label{Be3}
v_1(x_0,b)=v_2(x_0,a)=0\quad \forall x_0\in [0,T].
\end{equation}
We make the following  standard assumptions:
\begin{equation}\label{Deq6}
\tilde\rho_1,\tilde a_1, \tilde b_1\in C^1(\bar Q_+), \quad \tilde\rho_2,\tilde a_2,\tilde  b_2\in C^1(\bar Q_-),\quad  \tilde c_1\in  L^\infty(0,T; L^2(\Omega_+)),\quad \tilde c_2\in L^\infty(0,T; L^2(\Omega_-)),
\end{equation}
there exist a positive constant $\alpha_0$ such that
\begin{equation}\label{Deq7}
\tilde\rho_1(x)\ge \alpha_0>0,\quad
\tilde a_1(x)\ge \alpha_0>0 \quad \forall x\in Q_+,\quad \tilde \rho_2(x)\ge \alpha_0>0,\quad \tilde a_2(x)\ge \alpha_0>0 \quad \forall x\in\, Q_-.
\end{equation}
We have
\begin{proposition}\label {mika} Let $f_1\in L^2(Q_+), f_2\in L^2(Q_-), \tilde r\in L^2(0,T)$ and (\ref{Deq6}), (\ref{Deq7}) holds true. Then there exists function $\eta(x_1)\in C^2(\bar \Omega)$ , $\eta(x_1) <0$ on $\bar\Omega$ and  a constant $C_1$ independent of  $v=(v_1,v_2)$ such that
\begin{eqnarray}\label{mika1}
\sum_{\vert\alpha\vert\le 1}\Vert ((T-x_0)^{-3})^\frac {(3-2\vert\alpha\vert)}{2} \partial^\alpha  v_2\,
e^{\psi^*}\Vert_{L^2(Q_-)}  +\sum_{\vert\alpha\vert\le 1}\Vert ( (T-x_0)^{-3})^{\frac{5-2\vert\alpha\vert}{2}} \partial^\alpha  v_1\,
e^{\psi^*}\Vert_{L^2(Q_+)}\nonumber\\
+ \Vert ((T-x_0)^{-3})^\frac 32\partial_{x_1}^+v e^{\psi^*}\Vert_{L^2(0,T)}+\Vert ((T-x_0)^{-3})^\frac 12\partial_{x_1}^-v e^{\psi^*}\Vert_{L^2(0,T)}\nonumber\\+\Vert ((T-x_0)^{-3})^\frac 12\partial_{x_0} v e^{\psi^*}\Vert_{L^2(0,T)}+\Vert ((T-x_0)^{-3})^\frac 52 v e^{\psi^*}\Vert_{L^2(0,T)}\nonumber\\
\le C_1(\Vert (T-x_0)^{-3}f_1 e^{\psi^*}\Vert_{L^2(Q_+)}+\Vert (T-x_0)^{-3}f_2 e^{\psi^*}\Vert_{L^2(Q_-)}\nonumber\\+\Vert (T-x_0)^{-\frac 32}\tilde re^{\psi^*})\Vert_{L^2(0,T)}
+\Vert ((T-x_0)^{-3})^\frac 32\partial_{x_1} v e^{\psi^*}(\cdot,b)\Vert_{L^2(0,T)}),
\end{eqnarray} where $\psi^*(x)=\eta(x_1)/(T-x_0)^3.$
\end{proposition}

{\bf Proof.}
Making the change of variables $x_0\rightarrow T-x_0$ and setting $w_j(x)=v_j(T-x_0,x_1)$ we have
\begin{equation}\label{Ceq1}
 \rho^*_1\partial_{x_0} w_1- \partial_{x_1}^2 w_1+ b^*_1\partial_{x_1} w_1+ c^*_1 w_1= f^*_1\quad \mbox{in}\,\, Q_+,
\end{equation}
\begin{equation}
\label{Ce2}
\rho^*_2\partial_{x_0} w_2- \partial_{x_1}^2 w_2+b^*_2\partial_{x_1} w_2+c^*_2 w_2=f^*_2\quad \mbox{in}\,\, Q_-,
\end{equation}
\begin{equation}\label{Ce3}
w_1(x_0,0)-w_2(x_0,0)=(\partial_{x_1}w_1-\partial_{x_1}w_2-M\partial_{x_0} w_1)(x_0,0)-  r^*(x_0)=0 \quad \mbox{on} \,[0,T],
\end{equation}
\begin{equation}\label{CCe3}
w_1(x_0,b)=w_2(x_0,a)=0\quad \mbox{on} \,[0,T],
\end{equation}
where  $\rho^*_j(x)=\tilde \rho_j(T-x_0,x_1)/a_j(T-x_0,x_1), b^*_j(x)=\tilde b_j(T-x_0,x_1)/a_j(T-x_0,x_1), c^*_j(x)=\tilde c_j(T-x_0,x_1)/a_j(T-x_0,x_1), f^*_j(x)=\tilde f_j(T-x_0,x_1)/a_j(T-x_0,x_1), r^*(x)=\tilde r(T-x_0,x_1) $ and $j\in\{1,2\}.$

Our next step is to construct  of variables in domain $Q_+$ such that the equation (\ref{Ceq1}) keeps  the same form after change of variables but the new coefficient $\rho^*_1$ satisfies
$$
\rho^*_1=\rho^*_2\quad\mbox{on}\,\, [0,T]\times \{0\}.
$$ Let $F(x):C^{1,2}(\bar Q_+,\bar Q_+)$ be the diffeomorphism of $Q_+$ on $Q_+$ such that $F=(F_1,F_2)$ and
$F_1(x)=x_0$ on $Q_+$. In order to construct the function $F_2$ consider a function $q(x_0)\in C^1[0,T], q(x_0)>C>0$ on  $[0,T].$ Set $\kappa_0=\frac {b}{\Vert q\Vert_{C^0[0,T]}+40}.$  Let $\eta_1(x_1)\in C^\infty[0,b], $ $\eta_1(x_1)=x_1$ on $[0,\kappa_0]$, $\frac{d\eta_1}{dx_1}\ge 0$ on $[0,b]$, $\eta_1=b/5$ for $x_1\in  [\frac{4b}{5},b]$ and  $\frac{d\eta_1}{dx_1}> 0$ on $[\kappa_0,\frac {b}{\Vert q\Vert_{C^0[0,T]}+20}].$
Let
$$
\eta_2(x_1)=\left \{ \begin{matrix} 0\quad \mbox{for}\quad x_1\in [0,\frac {b}{\Vert q\Vert_{C^0[0,T]}+20}],\\\frac{b(x_1-b_1)^3}{(b-b_1)^3}\quad \mbox{for} \quad x_1\in [b_1,b].\end{matrix}\right.\quad  F_2(x)=q(x_0)\eta_1(x_1)+\eta_2(x).
$$
Then on $[0,T]\times [0,\tilde b]$ we have
$$
F^{-1}(\tilde x)=(\tilde x_0, \tilde x_1/q(\tilde x_0))\quad \mbox{and} \quad DF(x_0,0)=\left (\begin{matrix} 1&0\\ 0&q(x_0)  \end{matrix}\right).
$$
Denote $\tilde x=F(x)$ and $\tilde w_1( \tilde x)=w_1(F^{-1}(\tilde x)).$  Then $w_1(x)=\tilde w_1 (F(x))$. Therefore  on $Q_+$
$$
\partial_{x_0}w_1=\partial_{\tilde x_0}\tilde w_1 \partial_{x_0}F_1(x) +\partial_{\tilde x_1}\tilde w_1 \partial_{x_0}F_2(x)=\partial_{\tilde x_0}\tilde w_1 +\partial_{\tilde x_1}\tilde w_1 \partial_{x_0}F_2(x).
$$
In particular
$$
\partial_{x_0}w_1=\partial_{\tilde x_0}\tilde w_1 \quad \mbox{on}\quad [0,T]\times \{0\}.
$$

On $Q_+$ we have
$$
\partial_{x_1}w_1=\partial_{\tilde x_0}\tilde w_1 \partial_{x_1}F_1(x) +\partial_{\tilde x_1}\tilde w_1 \partial_{x_1}F_2(x)=\partial_{\tilde x_1}\tilde w_1 \partial_{x_1}F_2(x)
$$
and

$$
\partial^2_{x_1}w_1=\partial_{\tilde x_1}\tilde w_1 \partial^2_{x_1}F_2(x)+\partial^2_{\tilde x_0\tilde x_1}\tilde w_1\partial_{x_1}F_1(x) \partial_{x_1}F_2(x)+\partial^2_{\tilde x_1}\tilde w_1 (\partial_{x_1}F_2(x))^2=
$$
$$
\partial_{\tilde x_1}\tilde w_1 \partial^2_{x_1}F_2(x)+\partial^2_{\tilde x_1}\tilde w_1 (\partial_{x_1}F_2(x))^2.
$$
Therefore function $\tilde w_1$ satisfies the parabolic equation
$$
 \rho_0\partial_{\tilde x_0} \tilde w_1-\alpha(\tilde x)\partial^2_{\tilde x_1}\tilde w_1+\beta(\tilde x)\partial_{\tilde x_1}\tilde w_1+\tilde c\tilde w_1=f_0,
$$
where $ f_0=f^*_1\circ F^{-1}, \rho_0=\rho^*_1\circ F^{-1}$  and $\beta(\tilde x)=(b_1^*\partial_{x_1}F_2+\partial_{x_0}F_2-\partial^2_{x_1}F_2(x) )\circ F^{-1},\alpha(\tilde x)=(\partial_{x_1}F_1)^2\circ F^{-1},\tilde c=c_1^*\circ F^{-1}.$ After division of the new equation by $\alpha$ we have
 $$
 \rho^{*}_1\partial_{\tilde x_0} \tilde w_1-\partial^2_{\tilde x_1}\tilde w_1+b_1^*(\tilde x)\partial_{\tilde x_1}\tilde w_1+c_1^*\tilde w_1=f_1^*
$$
with $\rho^*_1=\rho_0/\alpha,$ $b^*_1=\beta/\alpha, c_1^*=\tilde c/\alpha$  and $f_1^*=f_0/\alpha.$
Observe that on $[0,T]\times \{0\}$
$$
\partial^2_{x_1}w_1=\partial^2_{\tilde x_1}\tilde w_1 q^2(x_0).`
$$
So
$$
\rho^*_1(\tilde x_0)=\rho^*_1(\tilde x_0)/q^2(\tilde x_0).
$$
Then taking $q^2(\tilde x_0)=\rho^*_2(x_0,0)/\rho^*_2(x_0,0)
$ we obtain that the function $\rho$ given by formula
$$
\rho(x)=\left \{ \begin{matrix} \rho^*_1\quad \mbox{for}\quad x\in Q_+,\\ \rho^*_2\quad \mbox{for} \quad x\in Q_-\end{matrix}\right.
$$ is continuous on $\bar Q.$
Equations (\ref{Ce3}) are transformed to
\begin{equation}\label{Ze3}
\tilde w_1( x_0,0)-w_2(\tilde x_0,0)=(q(x_0)\partial_{x_1}\tilde w_1-\partial_{x_1}w_2-M\partial_{x_0} \tilde w_1)(x_0,0)- r^*(x_0)=0 \quad \mbox{on}\, \,[0,T],
\end{equation}
Hence instead of proving the observability estimate to  system (\ref{Beq1}) - (\ref{Be3}) it suffices to prove the observability estimate for the following system:
\begin{equation}\label{gora1A}
P(x,D)u= \rho\partial_{x_0}u
- \partial^2_{x_1}u
+ b(x)\partial_{x_1} u +c(x)u=f
\quad\mbox{in}\,\,Q\setminus [0,T]\times\{0\},
\end{equation}
\begin{equation}\label{gora2A}
u(\cdot,a)=u(\cdot,b)=0\quad \mbox{on} \,[0,T],
\end{equation}

\begin{equation}\label{gora2A}
[u](x_0,\cdot)= -\partial^-_{x_1}u(x_0,0)+\mu(x_0)\partial^+_{x_1}u(x_0,0)-M\partial_{x_0}u(x_0,0)-r =0\quad \mbox{on} \quad [0,T],
\end{equation}
where
$$
u(x)=\left \{ \begin{matrix} \tilde w_1\quad \mbox{for}\quad x\in Q_+,\\ w_2\quad \mbox{for} \quad x\in Q_-\end{matrix}\right.,\quad\,b(x)=\left \{ \begin{matrix} b^*_1\quad \mbox{for}\quad x\in Q_+,\\ b^*_2\quad \mbox{for} \quad x\in Q_- \end{matrix}\right.,
$$
$$
c(x)=\left \{ \begin{matrix} c^*_1\quad \mbox{for}\quad x\in Q_+,\\ c^*_2\quad \mbox{for} \quad x\in Q_- \end{matrix}\right. ,\quad\, f(x)=\left \{ \begin{matrix} f^*_1\quad \mbox{for}\quad x\in Q_+,\\ f^*_2\quad \mbox{for} \quad x\in Q_-.\end{matrix}\right.
$$
Therefore the coefficients of equation (\ref{gora1A}) have the following regularity:
\begin{equation}\label{mila}
\rho\in C^1(\bar Q_+)\cap C^1(\bar Q_-)\cap
C^0(\bar Q) ,\quad \rho(x)>\beta>0\quad \mbox{on}\,\, Q,\quad
\end{equation}
\begin{equation}\label{mila1} \mu\in C^1[0,T], \quad \mu(x_0)>\beta>0\quad \mbox
{and}\quad b\in L^\infty(Q),\, c\in L^\infty(0,T;L^2(\Omega)).
\end{equation}
We set \begin{equation}\label{gopnikz}
\tilde\varphi(x_0)=\left\{ \begin{matrix}\frac{1}{x_0^3}\,\,\mbox{for}\,\, x_0\in [0,\frac T4]\\
\frac{1}{(T-x_0)^3}\,\,\mbox{for}\,\, x_0\in [\frac {3T}4,T]\end{matrix}\right. , \quad \varphi_*(x)=\left\{\begin{matrix} \varphi_2\quad \mbox{on}\,\, Q_+, \\ \varphi_1\quad \mbox{on}\,\, Q_-,\end{matrix}\right.
\end{equation}
where
\begin{equation}\label{gopnik1}
\varphi_j(x)=\frac{e^{\lambda\psi_j(x_1)}-e^{10000\lambda c_0}}{x_0^3(T-x_0)^3}, \quad c_0=\mbox{max}\{b,-a\},\quad j\in\{1,2\},
\end{equation} where $\lambda$ is a large positive parameter, $\tilde \varphi\in C^2[\frac {T}{8}, \frac{8T}{9}],$ and strictly positive on $[\frac {T}{8}, \frac{8T}{9}]$
  and \begin{equation}\label{gopnik2}\psi_1(x_1)=(x_1+10+c_0)^2\quad\mbox{ and}\quad \psi_2(x_1)=\psi_1(0)e^{(x_1+10+c_0)^2-(10+c_0)^2}.
  \end{equation}
  By (\ref{gopnikz}) - (\ref{gopnik2}) the following is true:
  \begin{equation}\label{barbos}
  \varphi_2(x)>\varphi_1(x)\quad \mbox{on}\,\, Q_+,\quad \varphi_1(x)>\varphi_2(x)\quad \mbox{on}\, \, Q_-, \,\, \varphi_1(x_0,0)=\varphi_2(x_0,0)\quad \forall x_0\in [0,T].
  \end{equation}
We introduce the Hilbert  space
$$
\Vert f\Vert_Y=\root \of{\Vert f\Vert^2_{L^2(Q_-)}+\Vert s\tilde \varphi f\Vert^2_{L^2(Q_+)}},
$$
the  operator
$$
\mbox{\bf P}(x,D)u=(\rho\partial_{x_0}u
- \partial^2_{x_1}u,-\partial^-_{x_1}u(\cdot,0)+\mu\partial^+_{x_1}u(\cdot,0)-M\partial_{x_0}u(\cdot,0)): \mathcal X\rightarrow Y\times L^2[0,T]
$$
and the Banach space
$$\mathcal X=\{u\vert u\in H^{1,2}(Q_+)\cap  H^{1,2}(Q_-), \mbox{\bf P}(x,D)u\in (L^2(Q_+)\cap L^2(Q_-))\times L^2(0,T),
$$
$$[u](\cdot,0)=0,\quad u(\cdot,0)\in H^1_0(0,T), \,u(\cdot,a)=u(\cdot,b)=0.\}$$

Denote $\mathcal B v=(\partial^+_{x_1} v,\partial^-_{x_1} v, v)(\cdot,0), $  and   $\mathcal Z(0,T)=L^2_{(s\tilde\varphi)^3}(0,T)\times L^2_{s\tilde\varphi}(0,T)\times H^{1,\tilde s}(0,T)\cap L^2_{(s\tilde \varphi)^5}(0,T).$
We have
\begin{proposition}\label{zoopa} Let $u\in \mathcal X$ and coefficients $\rho,\mu, b, c$ satisfy (\ref{mila}), (\ref{mila1}) and parameter $\lambda$ fixed sufficiently large.
There exists  $s_0>1$  and positive constant $C_2$ such that
 for all
$s\ge s_0$ the following estimate  holds true
\begin{eqnarray}\label{main1}
\sum_{\vert\alpha\vert\le 1}\Vert (s\widetilde\varphi)^\frac {(3-2\vert\alpha\vert)}{2} \partial^\alpha  u\,
e^{s\varphi_*}\Vert_{L^2(Q_-)}  +\sum_{\vert\alpha\vert\le 1}\Vert (s\widetilde\varphi)^{\frac{5-2\vert\alpha\vert}{2}} \partial^\alpha  u\,
e^{s\varphi_*}\Vert_{L^2(Q_+)}
+ \Vert \mathcal B (ue^{s\varphi_*})\Vert_{\mathcal Z(0,T)}\nonumber\\
\le C_2(\Vert (f e^{s\varphi_*}, re^{s\varphi_*})\Vert_{Y\times L^2_{s\tilde \varphi}(0,T)}
+\Vert (s\tilde\varphi)^\frac 32\partial_{x_1} u e^{s\varphi_*}(\cdot,b)\Vert_{L^2(0,T)}),
\end{eqnarray}
where $C_2$ is independent of $s.$
\end{proposition}


{\bf Proof.} Without loss of generality using the standard arguments (see e.g. \cite{H})  we can prove  an estimate (\ref{main1}) under assumption $b=c\equiv 0.$
First, by an argument based on the partition of unity (e.g., Lemma 8.3.1
in \cite{H}), it suffices to prove the inequality (\ref{main1})
locally, by assuming that
\begin{equation}\label{3.55}
\text{supp}\, u\subset B (x^*,\delta),
\end{equation}
where $B(x^*,\delta)$ is the ball in $\Bbb R^2$ of the radius $\delta>0$
centered at some point $x^*$.

Let $\tilde \theta\in C^\infty_0(\frac 12,2)$ be a nonnegative function such that
\begin{equation}\label{knight}\sum_{\ell=-\infty}^\infty\tilde \theta(2^{-\ell}t)=1
\quad  \mbox{for all}\,\, t\in \Bbb R^1.
\end{equation} (For existence of such a function $\tilde\theta$ see e.g. \cite{Sog}.)

Set $u_\ell(x)= u(x)\kappa_\ell(x_0)$ where
\begin{equation}\label{knight0}
\kappa_\ell(x_0)=
\tilde\theta \left (2^{-\ell}2^\frac{1}{ \theta(x_0)^\frac 14}\right ),
\end{equation} where
\begin{eqnarray}\label{Aziza}\theta\in C^\infty[0,T],\quad \theta\vert_{[0,T/4]}=x_0,\quad \theta\vert_{[3T/4,T]}=T-x_0,\nonumber\\ \partial_{x_0}\theta (x_0)<0 \,\,\mbox{on}\,\, (0,\frac{T}{2}), \quad \partial_{x_0}\theta (x_0)>0 \,\,\mbox{on}\,\, (\frac T2,T),\quad \partial^2_{x_0}\theta(\frac{T}{2})<0 .
\end{eqnarray}
Observe that it
suffices to prove the Carleman estimate (\ref{main1}) for the function
$u_\ell$ instead of $ u$ provided that the constant
$C_1$  and the function $s_0$ are independent of $\ell.$
Observe that if $G\subset \Bbb R^m$ is a bounded domain and $q\in L^2(G)$,
then there exist an independent constants $C_3$ and $C_4$
(see e.g. \cite{Sog}) such that
\begin{equation}\label{gorokn}
C_3\sum_{\ell=-\infty}^{\infty}\Vert \kappa_\ell q\Vert^2_{L^2(G)}
\le \Vert q\Vert^2_{L^2(G)}\le C_4\sum_{\ell=-\infty}^{\infty}\Vert \kappa_
\ell q\Vert^2_{L^2(G)}.
\end{equation}
Denote the norm on the left-hand side of (\ref{main1}) as $\Vert\cdot\Vert_*.$ Suppose that the estimate  (\ref{main1}) is true for any function $u_\ell$ with constants $C_1$ and $s_0$ independent of $\ell.$
By (\ref{gorokn})  for some constant $C_5$ independent of $s$ we have

\begin{equation}\label{zombi}
\Vert { u} e^{s\varphi_*}\Vert_*
=\Vert \sum_{\ell=-\infty}^{+\infty} { u}_\ell e^{s\varphi_*}\Vert_*\le
\sum_{\ell=-\infty}^{+\infty}\Vert  { u}_\ell e^{s\varphi_*}\Vert_*\le
C_5\sum_{\ell=-\infty}^{\infty}(\Vert \kappa_\ell \mbox{\bf P}(x,D)
{ u}
e^{s\varphi_*}\Vert^2_{Y\times L^2_{s\tilde \varphi}(0,T)}
\end{equation}
$$+\Vert e^{s\varphi_*}[\kappa_\ell,\mbox{\bf P}(x,D)]
{ u}
\Vert^2_{Y\times L^2_{s\tilde \varphi}(0,T)}
+\Vert( s\tilde\varphi)^\frac 32 \kappa_\ell\partial_{x_1} ue^{s\varphi_*}(\cdot,b)\Vert^2_{L^2(0,T)})^\frac 12.
$$
By (\ref{gorokn}) we obtain from (\ref{zombi}):
\begin{equation}\label{ozon}
\Vert  u e^{s\varphi_*}\Vert_*
\le
C_6(\Vert \mbox{\bf P}(x,D)
{u}
e^{s\varphi_*}\Vert^2_{Y\times L^2_{s\tilde \varphi}(0,T)}
\end{equation}
$$+\sum_{\ell=-\infty}^{\infty}\Vert e^{s\varphi_*}[\kappa_\ell,\mbox{\bf P}(x,D)]
{u}
\Vert^2_{Y\times L^2_{s\tilde \varphi}(0,T)}
+ \Vert( s\tilde\varphi)^\frac 32 \partial_{x_1}{u} e^{s\varphi_*}(\cdot,b)\Vert^2_{L^2(0,T)})^\frac 12.
$$
Using (\ref{knight0}) and (\ref{Aziza}) we
estimate the norm of the commutator $ [\kappa_\ell,\mbox{\bf P}(x,D)]$ we obtain
\begin{eqnarray}\label{ooo}
\sum_{\ell=-\infty}^{\infty}\Vert e^{s\varphi_*} [\kappa_\ell,\mbox{\bf P}(x,D)]{u}
\Vert^2_{Y\times L^2_{s\tilde \varphi}(0,T)}\le C_{7}\sum_{\ell=-\infty}
^{\infty}(
\Vert \partial_{x_0}\kappa_\ell u(\cdot,0)
e^{s\varphi_*}\Vert^2_{L^2_{s\tilde \varphi}(0,T)}+\Vert \partial_{x_0}\kappa_\ell u
e^{s\varphi_*}\Vert^2_Y)\nonumber\\
\le C_{8}\sum_{\ell=-\infty}
^{\infty}(
\Vert \widetilde \varphi^\frac {5}{12} \chi_{\mbox{supp}\, \kappa_\ell}{u}
e^{s\varphi_*}\Vert^2_{Y} +\Vert \widetilde \varphi^\frac {5}{12} \chi_{\mbox{supp}\, \kappa_\ell} u(\cdot,0)
e^{s\varphi_*}\Vert^2_{L^2_{s\tilde \varphi}(0,T)})\nonumber\\
\le C_{9}(\Vert \widetilde \varphi^\frac {5}{12} { u}
e^{s\varphi_*}\Vert^2_{Y}+\Vert \widetilde \varphi^\frac {5}{12} u(\cdot,0)
e^{s\varphi_*}\Vert^2_{L^2_{s\tilde \varphi}(0,T)})
.
\end{eqnarray}

From  (\ref{ozon}), and (\ref{ooo}) we obtain (\ref {main1}).

Now, without loss of generality we assume that
\begin{equation}\label{3.55}
\text{supp}\, { u}\subset B (x^*,\delta/2)\cap
\mbox{supp}\, \kappa_\ell,
\end{equation}
where $B(x^*,\delta)$ is the ball of the radius $\delta>0$
centered at some point $x^*.$ If $x^*$ does not belong to the set $[0,T]\times \{0\}$ the estimate (\ref{main1}) is proved in \cite{Im}.
More specifically if $\mbox{supp}\,  u \cap ([0,T]\times\{0\})=\emptyset$ there exists a constant $C_{10}=C_{10}(\delta, x^*)$ and $s_0=s_0(\delta, x^*)$ such that
\begin{eqnarray}\label{Lider}
\sum_{\vert\alpha\vert\le 1}\Vert (s\widetilde\varphi)^\frac {(3-2\vert\alpha\vert)}{2} \partial^\alpha  u\,
e^{s\varphi_*}\Vert_{L^2(Q_-)}  +\sum_{\vert\alpha\vert\le 1}\Vert (s\widetilde\varphi)^{\frac {5-\vert\alpha\vert}{2}} \partial^\alpha  u\,
e^{s\varphi_*}\Vert_{L^2(Q_+)}\\
\le C_{10}(\Vert (P(x,D)u)e^{s\varphi_*}\Vert_{Y}
+\Vert (s\tilde\varphi)^\frac 32\partial_{x_1} u e^{s\varphi_*}(\cdot,b)\Vert_{L^2(0,T)}),\nonumber
\end{eqnarray}
where $C_{10}$ is independent of $s.$
Therefore we have to consider the case
\begin{equation}\label{koma}
x^*=(x^*_0,0),\quad \text{supp}\, { u}\subset B (x^*,\delta/2)\cap
\mbox{supp}\, \kappa_\ell, \quad  B (x^*,\delta)\cap ([0,T]\times \{a,b\})
=\emptyset.
\end{equation}

For any function  $\varphi\in \{\varphi_1,\varphi_2\}$ we introduce the operator
$$
P_\varphi (x,D, \tilde s)=i\rho(x) D_0+(D_1+\vert \tilde s\vert i\phi_0(x,x^*))^2,\quad \phi_0(x,x^*)=\partial_{x_1}\varphi(x)/\tilde\varphi(x^*),\quad
\tilde s=s \tilde\varphi(x^*).
$$

For any $\xi_0\in \Bbb R^1\setminus\{0\} $ and $x\in Q$  we choose $\root \of {i\rho(x)\xi_0 }$ such that
\begin{equation}\label{xyi}
\mbox{Im}\, \root \of {i\rho(x)\xi_0} >0.
\end{equation} By (\ref{mila}) such choice is possible.
We define symbol $p_{\varphi}(x,\xi,\widetilde s)$ by formula
\begin{eqnarray}\label{3.34}
p_{\varphi}(x,\xi,\widetilde s) = i\rho(x)\xi_0+ (\xi_1+{i}\vert \tilde s\vert\phi_0(x,x^*))^2,
\end{eqnarray}

The zeros of the polynomial $p_{\varphi}(x,\xi,\widetilde s)$ with respect to variable $\xi_1$  for
$M(\xi_0, \widetilde s) \ge 1, $ and $x\in B (x^*,\delta)\cap
\mbox{supp}\, \kappa_\ell$ are
\begin{equation}\label{3.35}
\Gamma^\pm_\varphi(x,\xi_0,\widetilde s)
= (-{i}\vert \widetilde s\vert\widetilde\mu_\ell\varphi_{0}\kappa(\xi_0,\widetilde s)
+\alpha^\pm (x,
\xi_0,\widetilde s)),
\end{equation}
where
\begin{equation}\label{bobo}\widetilde\mu_\ell(x)=\eta_*(x)\sum_{k=\ell-20}^{\ell+20}\kappa_\ell(x_0),\quad
\eta_*\in C_0^\infty (B(x^*,2\delta)),\quad \eta_*\vert_{B(x^*,\delta)}=1,\end{equation}
the function $\kappa_\ell$ is given by (\ref{knight0}),
\begin{equation}\label{3.36}
\alpha^\pm(x,\xi_0,\widetilde s)=\pm\widetilde \mu_\ell(x)\kappa(\nu,\xi_0,\widetilde s)\root
\of{i\rho \xi_0}.                          \end{equation}
Next we construct the function $\kappa(\xi_0,\widetilde s)=\kappa(\nu,\xi_0,\widetilde s).$
Let
$\chi_\nu$ be a $C^\infty_0(\Bbb M)$ function
such that
$\chi_\nu$ is identically equal $1$ in some conic neighborhood  of the $( \xi^*_0,\tilde s^*)\in
\Bbb M$ and $\mbox{supp}\, \chi_\nu( \xi_0,\tilde s)\subset \mathcal O(\zeta^*,\delta_1).$
Assume that
\begin{equation}\label{book}
\kappa(\nu,\xi_0,\widetilde s)\vert_{\mbox{supp}\, \chi_\nu}=1,\quad \mbox{supp}\, \kappa(\nu,\xi_0,\widetilde s)
\subset \mathcal O(\zeta^*,2\delta_1), \,\,1\ge \kappa(\nu,\xi_0,\widetilde s)\ge 0\quad\mbox{on}\,\, \Bbb M.
\end{equation}

 We extend the function $\chi_\nu$ on $\Bbb R^{2}$ as follows :
$\chi_\nu(\xi_0/M^2(\xi_0, \widetilde s), \widetilde s/M(\xi_0, \widetilde s))$ for
$M(\xi_0,\widetilde s)>1$ and
$\chi_\nu(\xi_0/M^2(\xi_0, \widetilde s),\widetilde  s/M(\xi_0, \widetilde s))\kappa^*(
M(\xi_0,\widetilde s))$  for $M(\xi_0, \widetilde s)
<1$, where $\kappa^*(t)\in C^\infty(\Bbb R^1), \kappa^*(t)\ge 0, \kappa^*(t)=1$ for
$t\ge 1$ and $\kappa^*(t)=0$ for $t\in [0,1/2].$ In the similar way we extend the function $\kappa(\nu,\xi_0,\widetilde s)$ on $\Bbb R^{2}$. Denote by $\tilde \chi_\nu(x,D_0,\widetilde s)$ the pseudodifferential operator
with the symbol $\eta_\ell(x)\chi_\nu( \xi_0,\widetilde s)$ and
\begin{equation}\label{cloun}\eta_\ell(x)=\eta_{**}(x)
\sum_{k=-10}^{10}\kappa_{\ell+k}(x_0),\end{equation} where  $\eta_{**}\in C_0^\infty (B(x^*,\delta)),\quad \eta_{**}\vert_{B(x^*,\frac 34\delta)}=1.$
We set $$v_{\nu,\varphi}
=\tilde\chi_\nu(x, D_0,\widetilde s)v_\varphi\quad
\mbox{and} \quad  v_\varphi={ u}e^{s\varphi},$$
$$ r_{\nu,\varphi}
=\tilde\chi_\nu(x, D_0,\widetilde s)r_\varphi\quad
\mbox{and} \quad  r_\varphi={ r}e^{s\varphi}.$$

Let ${\mathcal O}$ be a domain in ${\Bbb R}^1.$

{\bf Definition.} {\it We say that the symbol
$a(x_0,\xi_0,\widetilde  s)\in C^{\widetilde k}(\overline{\mathcal O}
\times
{\Bbb R}^{2})$ belongs to the class
$C^{\widetilde k}_{cl}S^{\kappa,\widetilde  s}({\mathcal O})$ if

{\bf A}) There exists a compact set $K\subset\subset{\mathcal O}$
such that $a(x_0,\xi_0,\widetilde s)\vert_{{\mathcal O}
\setminus K}=0;$

{\bf B}) For any $\beta=(\beta_0,\beta_{1})$ there exists a
constant $C_\beta$
$$
\left\Vert \partial^{\beta_0}_{\xi_0}
\partial^{\beta_{1}}_{\widetilde s}
a(\cdot,\xi_0,\widetilde s)\right\Vert_{C^{\widetilde k}(\overline{\mathcal O})}\le
C_\beta\left(\widetilde s^2+\vert \xi_0\vert\right)
^{\frac{\kappa-\vert\beta\vert}{2}}\quad,
$$
where $\vert \beta\vert=2\beta_0+\beta_1$ and
$M(\xi_0,\widetilde s) \ge 1$;

{\bf C}) For any $N\in \Bbb N$ the symbol $a(x_0, \xi_0,\widetilde s)$ can be represented
as
$$
a(x_0, \xi_0,\widetilde s)=\sum_{j=1}^Na_j(x_0,
 \xi_0,\widetilde s)
+ R_N(x_0,\xi_0,\widetilde s)
$$
where the functions $a_j$ have the following properties: for any $\tilde\lambda>1$ and for all $(x_0,\xi_0,\widetilde  s)\in\{(x_0,
\xi_0,\widetilde  s)
\vert x_0\in K, \,M (\xi_0,\widetilde s) >1\}$
$$
a_j(x_0,{\tilde\lambda}^2 \xi_0,{\tilde\lambda}
\widetilde s)={\tilde\lambda}^{\kappa-j}a_j(x_0, \xi_0,\widetilde  s);
$$
for any multiindex $\beta$ and any and $(\xi_0,\widetilde  s)$ satisfying
$M(\xi_0,\widetilde s)\ge 1$ there exist a constant $C_\beta$ such that
$$
\left\Vert  \partial^{\beta_0}_{\xi_0}
\partial^{\beta_{1}}_{\widetilde s}
a_j(\cdot, \xi_0,\widetilde s)\right\Vert_{C^{\widetilde k}(\bar{\mathcal O})}\le
C_\beta\left(\widetilde s^2+\vert \xi_0\vert\right)
^{\frac{\kappa-j-\vert\beta\vert}{2}}
$$
where the term $R_N$ satisfies the estimate
$$
\Vert R_N(\cdot,\xi_0,\widetilde s)\Vert_{C^{\widetilde k}(\overline{\mathcal O})}\le
C_N(\widetilde s^2+\vert \xi_0\vert)^{\frac{\kappa-N}{2}}\quad
\forall (\xi_0,\widetilde  s)\,\, \mbox{satisfying}\,\, M
(\xi_0,\widetilde s)
 \ge 1.
$$}

For the symbol $a$, we introduce the semi-norm
\begin{eqnarray}
\pi_{C^{\widetilde k}({\mathcal O})}(a)=\sum_{j=1}^{\widehat N}\sup_{\vert \beta\vert\le \widehat N}
\sup_{\vert( \xi_0,\tilde s)\vert\ge 1}\left\Vert \frac{\partial^{\beta_0}}
{\partial\xi_0^{\beta_0}}
\frac{\partial^{\beta_{1}}}{\partial \widetilde s^{\beta_{1}}}
a_j(\cdot,\xi_0,\tilde s)\right\Vert_{C^{\widetilde k}(\overline{\mathcal
O})}/(1+\vert(\xi_0,\tilde s)\vert)^{\kappa-j-\vert\beta\vert} \nonumber \\
+\sup_{\vert (\xi_0,\tilde s)\vert\le 1}\Vert
a(\cdot,\xi_0,\tilde s)\Vert_{C^{\widetilde k}(\overline{\mathcal O})}. \nonumber
\end{eqnarray}
Obviously for any  $\widetilde k\in \{0,1\}$
\begin{equation}\label{zefir}
\pi_{C^{\widetilde k}({B(0,\delta(x^*))})}(\chi_\nu)\le C_{11}\widetilde\varphi^\frac {5\widetilde k}{12}(x^*).
\end{equation}
Obviously the  pseudodifferential operators with the symbols
$\Gamma_\varphi^\pm$ belongs \\
to the class $C^{\widetilde k}_{cl}S^{1,s}(B(0,\delta(x^*)))$  for any $\widetilde k\in \{0,1\}$
and
\begin{equation}\label{zigmund}
\pi_{C^{\widetilde k}(B(0,\delta(x^*)))}(\Gamma_\varphi^\pm)\le C_{12}\widetilde\varphi^\frac {5\widetilde k}{12}(x^*).
\end{equation}
By (\ref{zigmund}) and Lemma  8.1  of \cite{IM}
\begin{equation}\label{norm}
\Vert \Gamma_\varphi^\pm(\cdot,0, D_0,\tilde s)v_{\nu,\varphi}(\cdot,0)\Vert_{L^2(0,T)}\le C_{13}\Vert v_{\nu,\varphi}\Vert_{H^{\frac 12, \tilde s}(0,T)}.
\end{equation}
In some cases, we can represent the operator $
P_{\varphi}(x,D,\widetilde s)$ as a product of two first order
pseudodifferential operators.
\begin{proposition}(see e.g. \cite{IM})\label{gorokx1}  Let $v\in  \mathcal X,$ $\mbox{supp}\, v\subset
B(x^*,\delta)\cap \mbox{supp}\,\kappa_\ell, x^* \in \mbox{supp}\,\kappa_\ell$  and $ P_\varphi(x,D,\widetilde s)
\chi_\nu  v \in L^2(Q_+)\cap L^2(Q_-).$ Assume that  $\xi_0^*\ne 0$
and $ supp\, {\chi_\nu}\subset \mathcal O(\zeta^*,\delta_1).$
Then we can
factorize the operator $P_{\varphi}(x,D,\widetilde s)$ into the product of two
first order pseudodifferential operators:
\begin{eqnarray} \label{min}
P_{\varphi}(x,D,\widetilde s)v_{\nu,\varphi}=(D_1-\Gamma^-_\varphi
(x,{ D}_0,\widetilde s))
(D_1-\Gamma^+_\varphi(x,{ D}_0,\widetilde s))v_{\nu,\varphi}
+T_{+,\varphi} v_{\nu,\varphi}=\\(D_1-\Gamma^+_\varphi
(x,{ D}_0,\widetilde s))
(D_1-\Gamma^-_\varphi(x,{ D}_0,\widetilde s))v_{\nu,\varphi}+T_{-,\varphi} v_{\nu,\varphi},\nonumber
\end{eqnarray}
Operators  $T_{\pm,\varphi} : H^{\frac 12 ,1,\widetilde s}([0,b]\times {\Bbb R}^{1})
\rightarrow L^2(0,b;L^2({\Bbb R}^1))\cap H^{\frac 12 ,1,\widetilde s}([a,0]\times {\Bbb R}^{1})
\rightarrow L^2(a,0;L^2({\Bbb R}^1))$ satisfy estimates
\begin{equation}\label{lokom}
\Vert T_{\pm,\varphi} v_{\nu,\varphi}\Vert_{L^2(0,b;L^2({\Bbb
R}^1))}\le C_{14}\widetilde\varphi^\frac {5}{12}(x^*)\Vert v_{\nu,\varphi}\Vert_{H^{\frac 12,1,\tilde s}(Q_+)}.
\end{equation}
and
\begin{equation}\label{lokom}
\Vert T_{\pm,\varphi} v_{\nu,\varphi}\Vert_{L^2(a,0;L^2({\Bbb
R}^1))}\le C_{15}\widetilde\varphi^\frac {5}{12}(x^*)\Vert v_{\nu,\varphi}\Vert_{H^{\frac 12,1,\tilde s}(Q_-)}.
\end{equation}
\end{proposition}

Denote by
$
V^\pm_{\nu,\varphi}=(D_1-\Gamma^\pm_{\varphi}(x,{D}_0,
\widetilde s))
v_{\nu,\varphi }
$ the function with domain in $\bar Q_+$ and by $
U^\pm_{\nu,\varphi}=(D_1-\Gamma^\pm_{\varphi}(x,{D}_0,
\widetilde s))
v_{\nu,\varphi }
$ the function with domain in $\bar Q_-.$

Let us consider the initial value problems
\begin{equation}\label{barmalei0}
(D_1-\Gamma^-_\varphi(x, D_0,\widetilde s))V^+_{\nu,\varphi}=-T_{+,\varphi}v_{\nu,\varphi}+P_{\varphi}(x,D,\widetilde s)v_{\nu,\varphi} \quad x\in [0,b]\times \Bbb R^1.
\end{equation}
and
\begin{equation}\label{barmalei01}
(D_1-\Gamma^-_\varphi(x, D_0,\widetilde s))V^+_{\nu,\varphi}=-T_{+,\varphi}v_{\nu,\varphi}+P_{\varphi}(x,D,\widetilde s)v_{\nu,\varphi} \quad x\in [0,b]\times \Bbb R^1.
\end{equation}
For solutions of these problems, we can prove an a priori estimate.

\begin{proposition}\label{gorokx}(see e.g. \cite{IM}) Let $\xi_0^*\ne 0,\mbox{supp}\, \chi_\nu\in \mathcal O(\zeta^*,\delta_1).$ There exists a
constant $C_{16}>0$ such that
\begin{eqnarray} \label{min1}
\Vert V^+_{\nu,\varphi}(\cdot, 0)\Vert_{H^{\frac{1}{4}}(0,T)\cap L^2_{ \tilde s}(0,T)}+\Vert V^+_{\nu,\varphi}\Vert_{H^{\frac{1}{2},1,\tilde s}(Q_+)}\\\le
C_{16}(\tilde\varphi^\frac{5}{12}(x^*)\Vert  v\Vert_{H^{\frac 12,1,\tilde s}(Q_+) }+\Vert P_{\varphi}(x,D,\widetilde s)v_{\nu,\varphi}\Vert_{L^2(Q_+)}).\nonumber
\end{eqnarray}
and
\begin{eqnarray} \label{min11}
\Vert V^+_{\nu,\varphi}\Vert_{H^{\frac{1}{2},1,\tilde s}(Q_-)}\le C_{17}(\tilde\varphi^\frac{5}{12}(x^*)\Vert  v\Vert_{H^{\frac 12,1,\tilde s}(Q_-) }\\
+\Vert P_{\varphi}(x,D,\widetilde s)v_{\nu,\varphi}\Vert_{L^2(Q_+)}+\Vert  V^+_{\nu,\varphi}(\cdot, 0)\Vert_{H^{\frac{1}{4},\tilde s}[0,T]\cap L^2_{ \tilde s}(0,T)}).\nonumber
\end{eqnarray}
\end{proposition}

Consider the initial value problem:
\begin{equation}
(D_1-r(x, D_0,\widetilde s))W=p \quad x\in [-\delta,0]\times \Bbb R^1,
\quad W\vert_{x_1=-\delta}=0.
\end{equation}
We have
\begin{proposition}\label{gorokxY}(see e.g. \cite{IM}) Let $\xi_0^*\ne 0,\mbox{supp}\, \chi_\nu\in \mathcal O(\zeta^*,\delta_1),  W \in H^{\frac 12,1}(\Bbb R^1\times [-\delta,0]),p\in L^2(\Bbb R^1\times [-\delta,0]).$ Let  for each $x_1\in [-\delta,0]$  symbol $r(x_0,x_1,\xi_0,\widetilde  s)\in C^1_{cl}S^{1,\widetilde  s}((0,T)
) $ for all $x_1\in (-\delta,0)$ and there exist a constant $C_{18}$ such that
$$
r(x,\xi_0,\tilde s)\ge C_{18}M(\xi_0,\tilde s)\quad (x,\xi_0,\tilde s)\in (B(x^*,\frac 34 \delta)\cap \mbox{supp}\, \sum_{k=-8}^8 \kappa_{k+\ell})\times \Bbb R^2
$$
 and
 $$
 \pi_{C^1(0,T)
}  (r(\cdot,x_1,\cdot,\cdot))\le C_{19}\widetilde\varphi^\frac {5}{12}(x^*) \quad \forall x_1\in [-\delta,0].
$$
Then here exists a
constant $C_{20}>0$ such that
\begin{equation}
\Vert W(\cdot, 0)\Vert_{H^{\frac{1}{4}}[0,T]\cap L^2_{\vert \tilde s\vert}(0,T)}+\Vert W\Vert_{H^{\frac{1}{2},1,\tilde s}((-\delta,0)\times \Bbb R^1))}\\\le
C_{20}\Vert p\Vert_{L^2( (-\delta,0)\times \Bbb R^1)}.\nonumber
\end{equation}
\end{proposition}

Now we obtain couple subelliptic estimates for the operator $P_\varphi(x,D,\widetilde s)$ on domains $Q_\pm.$

\begin{proposition}(see e.g. \cite{IM})\label{opana} Let parameter $\lambda$ be large enough and fixed,  $w\in \mathcal X,$ $\mbox{supp}\, w\subset
B(x^*,\delta)\cap \mbox{supp}\,\kappa_\ell$  and $P_\varphi(x,D,\widetilde s)
\chi_\nu w\in L^2(Q_\pm).$ Then
there exist positive constants $\delta(x^*), C_{21}, C_{22}$
independent of $\tilde s$
 such that for all
and $\tilde s\ge s_0$ we have
\begin{eqnarray}\label{klop}
C_{21}\int_{ Q_+}(\vert\tilde  s\vert \vert \partial_{x_1}\tilde \chi_\nu w\vert^2
+ \vert\tilde  s\vert^3 \vert \tilde \chi_\nu w\vert^2)dx
- \mbox{Re}\,
\int_{\Bbb R^1}  \partial^+_{x_1}\chi_\nu w
\overline{\rho\partial_{x_0}\tilde  \chi_\nu w}\vert_{x_1=0} d x_0\\+\int_{\Bbb R^1}
(\vert\widetilde  s\vert\phi_{0}(x^*,x^*)\vert\partial^+_{x_1}
\tilde \chi_\nu w\vert^2 +\vert\widetilde
s\vert^3 \varphi^3_{0}(x^*,x^*)\vert\tilde \chi_\nu w\vert^2)\vert_{x_1=0}d x_0\nonumber\\
\le  \Vert P_\varphi(x,D,\widetilde s)\tilde \chi_\nu w\Vert^2_{L^2(Q_+)}
+ \epsilon_1(
\delta)
\Vert(\tilde s^\frac 12\partial^+_{x_1}(\tilde \chi_\nu w),\tilde s^\frac 32\tilde \chi_\nu w)(\cdot,0) \Vert^2_{L^2(\Bbb R^1)\times L^2(\Bbb R^1)},\nonumber
\end{eqnarray}
where $\epsilon_1(\delta)\rightarrow +0$ as $\delta\rightarrow
+0$
and
\begin{eqnarray}\label{klop1}
C_{21}\int_{ Q_-}(\vert \tilde s\vert\vert \partial_{x_1}\tilde \chi_\nu w\vert^2
+ \vert \tilde s\vert^3\vert\tilde  \chi_\nu w\vert^2)dx
\\
\le  \Vert P_\varphi(x,D,\widetilde s)\tilde \chi_\nu w\Vert^2_{L^2(Q_-)}+ C_{22}(
\vert\int_{\Bbb R^1}
 \partial^-_{x_1}\tilde \chi_\nu w
\overline{\rho\partial_{x_0}\tilde  \chi_\nu w}\vert_{x_1=0} d x_0\vert\nonumber\\+\int_{\Bbb R^1}
(\vert\widetilde  s\vert\varphi_{0}\vert\partial^-_{x_1}
(\tilde \chi_\nu w)\vert^2 +\vert\widetilde
s\vert^3 \varphi^3_{0}\vert\tilde \chi_\nu w\vert^2)\vert_{x_1=0}d x_0).
\nonumber
\end{eqnarray}

\end{proposition}

We use the following proposition proved in (\cite{IM}):
\begin{proposition}\label{gopnik} Let $-\infty<\alpha<\tilde a<\tilde b<\beta<+\infty,$ $p\in \Bbb N_+$ and  $\mbox{supp} \,\mbox{\bf v}\subset  [\tilde a, \tilde b].$
Then there exists an independent constant $C_{23}$ such that
\begin{equation}\label{zima1}\Vert M^p( D_0,\tilde s)\mbox{\bf v}
\Vert_{L^2 ([-R,R]\setminus [\alpha,\beta])}\le \frac{C_{23}}{(\min \{\tilde a-\alpha, \beta-\tilde b\})^p}\Vert \mbox{\bf v}\Vert
_{L^2(\Bbb R^1)}.
\end{equation}
\end{proposition}

We apply the  Proposition \ref{gopnik} in order to estimate  the  $H^{\frac 12,\widetilde s}\cap L^2_{\tilde s^2}$ norm of the function $(1-\eta_\ell)\chi_\nu( D_0, \widetilde s) v_\varphi(\cdot,0).$ Observe that by (\ref{Aziza}), (\ref{knight0}) for all sufficiently large  $\ell$
\begin{equation}\label{zima}
\mbox{supp}\, v_\varphi(\cdot,0)\subset [(\ell+2)^{-4}, (\ell-2)^{-4}]\cup [T-(\ell-2)^{-4}, T-(\ell+2)^{-4}]
\end{equation}
and
\begin{equation}\label{zima2}
\mbox{supp}\,(1- \eta_\ell(\cdot,0))\subset [0,T]\setminus[(\ell+11)^{-4}, (\ell-11)^{-4}]\cup [T-(\ell-11)^{-4}, T-(\ell+11)^{-4}].
\end{equation}
Therefore, by (\ref{zima}) and (\ref{zima2}) for all sufficiently  large $\ell$
\begin{eqnarray}\label{cosmos}
\Vert (1-\eta_\ell)\chi_\nu( D_0, \widetilde s) v_\varphi(\cdot,0)
\Vert_{H^{\frac 12,\widetilde s} (-R,R)\cap L^2_{\tilde s^2}(-R,R)}\\\le C_{24}\ell^{5}\Vert v_\varphi (\cdot,0)\Vert
_{L^2(\Bbb R^1)} \le C_{25}\widetilde\varphi^\frac{5}{12}(x^*)\Vert v_\varphi(\cdot,0)\Vert
_{L^2(\Bbb R^1)}.\nonumber
\end{eqnarray}
Here in the last inequality we used  the fact that
$\ell^5\le C_{26}\theta^\frac 54(x^*) \le C_{27}\varphi^\frac{5}{12}(x^*).
$
By arguments, same as in Lemma 8.5  of \cite{IM} we obtain
\begin{equation}\label{po1}
\vert \widetilde s\vert\Vert (1-\eta_\ell)\chi_\nu( D_0, \widetilde s)v_\varphi (\cdot,0)
\Vert_{H^{\frac 12,\widetilde s} (\Bbb R^1\setminus [-R,R])\cap L^2_{\tilde s^2}(\Bbb R^1\setminus [-R,R])}\le{ C_{28}}
\Vert v_\varphi(\cdot,0)\Vert_{L^2(\Bbb R^1)}.
\end{equation}
By (\ref{cosmos}) and (\ref{po1})
\begin{equation}\label{nikolas}\Vert (1-\eta_\ell)\chi_\nu( D_0, \widetilde s) v_\varphi (\cdot,0)
\Vert_{H^{\frac 12 ,\widetilde s} (\Bbb R^1)\cap L^2_{\tilde s^2}(\Bbb R^1)}\le C_{29}\widetilde \varphi^\frac{5}{12}(x^*)\Vert
v_\varphi(\cdot,0)\Vert_{L^2(\Bbb R^1)}.
\end{equation}

We introduce three sets
$$\mathcal Z_{\varphi,\pm}(\ell)=\{(x,\xi_0,\tilde s)\in \mbox{supp}\, \kappa_\ell\times \{0\} \times \Bbb M\vert \quad \pm( \vert\tilde s\vert\phi_0(x_0,0,x^*)-\mbox{Im}\, \root \of {i\rho(x_0,0)\xi_0}) >0 \}
$$
and
$$\mathcal Z_{\varphi,0}(\ell)=\{(x,\xi_0,\tilde s)\in \mbox{supp}\, \kappa_\ell\times \{0\} \times \Bbb M\vert\quad  \vert\tilde s\vert\phi_0(x_0,0,x^*)=\mbox{Im}\, \root \of {i\rho(x_0,0)\xi_0} =0 \}.
$$

If $\zeta^*\in \mathcal Z_{\varphi,0}(\ell)$ or $\zeta^*\in \mathcal Z_{\varphi,+}(\ell)$  and $\supp \,\chi_\nu \subset \mathcal O(\zeta^*,\delta_1)$ we have

\begin{eqnarray}\label{Zzara00}
\Vert \partial_{x_0} v_{\nu,\varphi}(\cdot,0)\Vert_{L^2(0,T)}\le \Vert \partial_{x_0}\eta_\ell v_{\nu,\varphi}(\cdot,0)\Vert_{L^2(0,T)}+\Vert\xi_0\chi_{\nu}Fv(\cdot,0)\Vert_{L^2(\Bbb R^1)}\\
\le C_{30}(\tilde\varphi^{\frac {5}{12}}(x^*)\Vert v_{\varphi}(\cdot,0)\Vert_{L^2(0,T)}+
\vert \tilde s\vert^2\Vert\chi_{\nu}Fv_\varphi(\cdot,0)\Vert_{L^2(\Bbb R^1)})\nonumber\\
\le C_{31}(\tilde\varphi^{\frac {5}{12}}(x^*)\Vert v_{\varphi}\Vert_{L^2(0,T)}+
\vert \tilde s\vert^2\Vert v_{\nu,\varphi}(\cdot,0)\Vert_{L^2(\Bbb R^1)}+\vert \tilde s\vert^2\Vert (1-\eta_\ell)\chi_\nu(D_0,\widetilde s) v_{\varphi}(\cdot,0)\Vert_{L^2(\Bbb R^1)})\nonumber\\
\le C_{32}(\vert\tilde s\vert\tilde\varphi^{\frac {5}{12}}(x^*)\Vert v_{\varphi}(\cdot,0)\Vert_{L^2(0,T)}+
\vert \tilde s\vert^2\Vert v_{\nu,\varphi}(\cdot,0)\Vert_{L^2(\Bbb R^1)}).\nonumber
\end{eqnarray}
Here to get the last inequality we used (\ref{nikolas}).

We will use the following proposition:
\begin{proposition} (\cite{IM})Let $\mbox{w}\in L^2(Q)$ and $\mbox{supp}\,\mbox{w} \subset\mbox{supp} \, \tilde\mu_\ell(x),$ where function $\tilde\mu_\ell$ defined by (\ref{bobo}).
Then there exists a constant $C_{33}$ such that
\begin{equation}\label{nikolasZ}\Vert (1-\eta_\ell)\chi_\nu( D, \widetilde s) \mbox{w}
\Vert_{H^{\frac 12 ,1,\widetilde s} (\Bbb R^2)}\le C_{33}\widetilde \varphi^\frac{5}{12}(x^*)\Vert
\mbox{w}\Vert_{L^2(\Bbb R^2)}.
\end{equation}
\end{proposition}

If $\zeta^*\in \mathcal Z_{\varphi,0}(\ell)$ or $\zeta^*\in \mathcal Z_{\varphi,+}(\ell)$  and $\supp \,\chi_\nu \subset \mathcal O(\zeta^*,\delta_1)$ using  inequality (\ref{nikolasZ}) we have
\begin{eqnarray}\label{ladymikki}
\Vert v_{\nu,\varphi}\Vert_{H^{\frac 12,0}(Q_\pm)}\le \Vert\chi_\nu(D_0,\tilde s) v_{\varphi}\Vert_{H^{\frac 12,0}(Q\pm)}+\Vert (1-\eta_\ell) \chi_\nu(D_0,\tilde s) v_{\varphi}\Vert_{H^{\frac 12,0}(Q_\pm)}=\nonumber\\
 \le C_{34}(\Vert\root\of {\vert\xi_0\vert}\chi_\nu(\xi_0,\tilde s)F v_{\varphi}\Vert_{L^2(\Bbb R^1\times I_\pm)}+\Vert (1-\eta_\ell) \chi_\nu(D_0,\tilde s) v_{\nu,\varphi}\Vert_{H^{\frac 12,0}(Q\pm)})\nonumber\\
\le  C_{35}\Vert\tilde s\chi_\nu(\xi_0,\tilde s)F v_{\varphi}\Vert_{L^2(\Bbb R^1\times I_\pm)}+\Vert (1-\eta_\ell)\chi_\nu(D_0,\tilde s)  v_{\varphi}\Vert_{H^{\frac 12,0}(Q_\pm)}\nonumber\\
 \le C_{36}(\vert\tilde s\vert\Vert  v_{\nu,\varphi}\Vert_{L^2(Q_\pm)}+\Vert (1-\eta_\ell) \chi_\nu(D_0,\tilde s) v_{\varphi}\Vert_{H^{\frac 12,0}(Q_\pm)\cap L^2_{\tilde s}(Q_\pm)})\nonumber\\
 \le C_{37}(\vert\tilde s\vert\Vert  v_{\nu,\varphi}\Vert_{L^2(Q_\pm)}+ \tilde\varphi^\frac{5}{12}(x^*)\Vert  v_{\varphi}\Vert_{ L^2(Q_\pm)}).
\end{eqnarray}
We consider three cases

{\bf Case 1.} Let $(x^*,\xi_0^*,\tilde s^*)\in \mathcal Z_{\varphi,+}(\ell)$ and $ \mbox{supp}\, \chi_\nu\subset \mathcal O(\zeta^*,\delta_1(\zeta^*)).$
We observe that by (\ref{mila}) and (\ref{xyi})
$$
\mbox{Im}\, \root \of {i\rho(x^*_0,0)\xi^*_0}=\frac{1}{\root\of{2}}\root\of{\rho(x^*_0,0)}\root\of{\vert\xi_0^*\vert}.
$$

Therefore  there exist positive  $\delta^0_1$ such that  for any $\delta_1\in (0,\delta^0)$ one can find positive $\delta_2\in (0,\frac 14)$ such that
\begin{equation}\label{monstr}
\rho (x)\vert\xi_0\vert < (2-\delta_2)\tilde s^2\phi^2_0(x,x^*)\quad\forall (x,\xi_0,\tilde s)\in B(x^*,\delta)\times \mathcal O (\zeta^*,\delta_1(\zeta^*)).
\end{equation}
We start from the estimate of the following boundary integral
\begin{eqnarray}\label{1}
\vert \mbox{Re}\,
\int_{\Bbb R^1} \partial^+_{x_1} v_{\nu,\varphi}
\overline{\rho\partial_{x_0}v_{\nu,\varphi}}\vert_{x_1=0} d x_0\vert\le \vert  \mbox{Re}\,\int_{\Bbb R^n} \partial^+_{x_1} v_{\nu,\varphi}
\overline{\rho(x^*)\partial_{x_0}v_{\nu,\varphi}}\vert_{x_1=0} d x_0\vert\nonumber\\
+\vert \mbox{Re}\, \int_{\Bbb R^1} \partial^+_{x_1} v_{\nu,\varphi}
\overline{(\rho(x^*)-\rho(x))\partial_{x_0}v_{\nu,\varphi}}\vert_{x_1=0} d x_0\vert=\mathcal I_1+\mathcal I_2.
\end{eqnarray}
We estimate integrals  $\mathcal I_j$ separately. Short computations, (\ref{monstr}), (\ref{knight0}), (\ref{cloun})  and (\ref{nikolas})  imply

\begin{eqnarray}\label{2}
\mathcal I_2=\sup_{x\in \mbox{supp} \,v_{\nu,\varphi}}\vert \rho(x^*)-\rho(x)\vert  \int_{\Bbb R^n} \vert \partial^+_{x_1} v_{\nu,\varphi}
\overline{\partial_{x_0}v_{\nu,\varphi}}\vert_{x_1=0}\vert d x_0\nonumber
\\
\le \sup_{x\in \mbox{supp} \,v_{\nu,\varphi}}\vert \rho(x^*)-\rho(x)\vert \Vert\partial^+_{x_1} v_{\nu,\varphi}(\cdot,0)\Vert_{L^2(0,T)}(\Vert \xi_0\chi_\nu(\xi_0,\tilde s) Fv_{\varphi}(\cdot,0)\Vert_{L^2(\Bbb R^1)}+\Vert \partial_{x_0}\eta_\ell \chi_\nu v_\varphi\Vert_{L^2(\Bbb R^1)})
\nonumber\\
\le C_{38}\sup_{x\in \mbox{supp} \,v_{\nu,\varphi}}\vert \rho(x^*)-\rho(x)\vert \Vert\partial^+_{x_1} v_{\nu,\varphi}(\cdot,0)\Vert_{L^2(0,T)}(\tilde s ^2\Vert v_{\nu,\varphi}(\cdot,0)\Vert_{L^2(0,T)}
\nonumber\\
+\tilde s ^2\Vert (1-\eta_\ell)\chi_\nu v_{\varphi}(\cdot,0)\Vert_{L^2(0,T)}+\Vert \partial_{x_0}\eta_\ell \chi_\nu v_\varphi\Vert_{L^2(\Bbb R^1)})
\nonumber\\
\le C_{39}\sup_{x\in \mbox{supp} \,v_{\nu,\varphi}}\vert \rho(x^*)-\rho(x)\vert(\int_{\Bbb R^1}
(\vert\widetilde  s\vert\phi_{0}(x^*,x^*)\vert\partial^+_{x_1}
v_{\nu,\varphi}\vert^2 +\vert\widetilde
s\vert^3 \phi^3_{0}(x^*,x^*)\vert v_{\nu,\varphi}\vert^2)\vert_{x_1=0}d x_0\nonumber \\+\vert\tilde s\vert\tilde \varphi^\frac 56 (x^*)\Vert v_\varphi(\cdot,0)\Vert^2_{L^2(0,T)}).\nonumber
\end{eqnarray}

Since function $\rho$ is continuous function taking parameter $\delta $ in (\ref{3.55}) sufficiently small we obtain
\begin{equation}\label{2}
\mathcal I_2=
\frac{\delta_2}{6}(\int_{\Bbb R^1}
(\vert\widetilde  s\vert\phi_{0}(x^*,x^*)\vert\partial^+_{x_1}
v_{\nu,\varphi}\vert^2 +\vert\widetilde
s\vert^3 \phi^3_{0}(x^*,x^*)\vert v_{\nu,\varphi}\vert^2)\vert_{x_1=0}d x_0+C_{40}\vert\tilde s\vert\tilde \varphi^\frac 56 (x^*)\Vert v_\varphi(\cdot,0)\Vert^2_{L^2(0,T)}).
\end{equation}

Next we estimate $\mathcal  I_1$
$$
\mathcal I_1\le  \vert \int_{\Bbb R^1} \partial^+_{x_1} v_{\nu,\varphi}
\overline{\rho(x^*)\partial_{x_0}v_{\nu,\varphi}}\vert_{x_1=0} d x_0\vert\le \vert \int_{\Bbb R^1} \vert F \partial^+_{x_1} v_{\nu,\varphi}
\overline{\rho(x^*) \xi_0 \chi_{\nu}(\xi_0,\tilde s)Fv_{\nu,\varphi}}\vert\vert_{x_1=0} d \xi_0\vert
$$
$$
+\vert \int_{\Bbb R^1} \partial^+_{x_1} v_{\nu,\varphi}
\overline{\rho(x^*)\partial_{x_0}\eta_\ell\chi_\nu (D_0,\tilde s) v_{\varphi}}\vert_{x_1=0} d x_0\vert .
$$
By (\ref{monstr})
\begin{eqnarray}\label{3}
\mathcal I_1\le ( 2-\delta_2)  \Vert\partial^+_{x_1} v_{\nu,\varphi}(\cdot,0)\Vert_{L^2(0,T)}
  \vert\tilde s\vert\phi_0(x^*,x^*)\Vert\chi_{\nu}(\xi_0,\tilde s)Fv_{\varphi}(\cdot,0)  \Vert_{L^2(\Bbb R^1)}\nonumber\\+\tilde \varphi^\frac {5}{12}(x^*)\Vert v_\varphi(\cdot,0)\Vert_{L^2(0,T)}\Vert \partial^+_{x_1} v_{\nu,\varphi}(\cdot,0)\Vert_{L^2(0,T)}
 \nonumber\\
 =  ( 2-\delta_2)   \Vert\partial^+_{x_1} v_{\nu,\varphi}(\cdot,0)\Vert_{L^2(0,T)}
 \vert\tilde s\vert\phi_0(x^*,x^*)\Vert v_{\nu,\varphi}(\cdot,0)  \Vert_{L^2(0,T)}+\nonumber
 \\
 ( 2-\delta_2)  \Vert\partial^+_{x_1} v_{\nu,\varphi}(\cdot,0)\Vert_{L^2(0,T)}\vert\tilde s\vert^2\phi^2_0(x^*,x^*)\Vert(1-\eta_\ell)\chi_{\nu}(D_0,\tilde s)v_{\varphi}\Vert_{L^2(\Bbb R^1)}\nonumber\\+C_{41}\tilde \varphi^\frac {5}{12}(x^*)\Vert v_\varphi(\cdot,0)\Vert_{L^2(0,T)}\Vert \partial^+_{x_1} v_{\nu,\varphi}(\cdot,0)\Vert_{L^2(0,T)}
 \nonumber\\
\le
(1-\delta_2/5)\int_0^T (\vert \tilde s\vert\phi_0(x^*,x^*)\vert\partial^+_{x_1}v_{\nu,\varphi}\vert^2+\vert\tilde s\vert^3\phi_0^3(x^*,x^*)\vert v_{\nu,\varphi}\vert^2)(x_0,0)dx_0\nonumber\\+C_{42}\vert\tilde s\vert \tilde \varphi^\frac{5}{6}(x^*)\Vert v_\varphi(\cdot,0)\Vert^2_{L^2(0,T)}.\nonumber
\end{eqnarray}
So by (\ref{1})-(\ref{3})
\begin{eqnarray}\label{groml}
\vert \mbox{Re}\,
\int_{\Bbb R^1} \partial_{x_1}^+ v_{\nu,\varphi}
\overline{\rho\partial_{x_0}v_{\nu,\varphi}}\vert_{x_1=0} d x_0\vert\\\le (1-\delta_2/5)\int_0^T (\vert \tilde s\vert\phi_0(x^*,x^*)\vert\partial^+_{x_1}v_{\nu,\varphi}\vert^2+\vert\tilde s\vert^3\phi_0^3(x^*,x^*)\vert v_{\nu,\varphi}\vert^2)(x_0,0)dx_0\nonumber\\+C_{43}\vert\tilde s\vert\tilde\varphi^\frac 56(x^*)\Vert v_\varphi(\cdot,0)\Vert^2_{L^2(0,T)}).\nonumber
\end{eqnarray}
By (\ref{groml}) and inequality (\ref{klop}) of Proposition \ref{opana}
\begin{eqnarray}\label{zara}
\int_0^T (\vert \tilde s\vert \vert\partial^+_{x_1}v_{\nu,\varphi}\vert^2+\vert\tilde s\vert^3\vert v_{\nu,\varphi}\vert^2)(x_0,0)dx_0+\vert\tilde s\vert\Vert v_{\nu,\varphi} \Vert^2_{H^{0,1,\tilde s}(Q_+)}\\\le C_{44}(\Vert P_{\varphi}(x,D,\widetilde s)v_{\nu,\varphi}\Vert^2_{L^2(Q_+)}+\vert\tilde s\vert\tilde\varphi^\frac 56(x^*)\Vert v_\varphi(\cdot,0)\Vert^2_{L^2(0,T)}).\nonumber
\end{eqnarray}
Equation
$$
\mathcal R u =-\partial^-_{x_1}u(x_0,0)+\mu\partial^+_{x_1}u (x_0,0)=M\partial_{x_0}u (x_0,0)+ r(x_0)
$$ in terms of functions $v_{\nu,\varphi}$ and  $r_{\nu,\varphi}$ one can write down on $[0,T]$ as
 \begin{eqnarray}\label{AAgora2A}
\mathcal R v_{\nu,\varphi} +[\chi_\nu,\mathcal R]v_{\varphi}=M\partial_{x_0}v_{\nu,\varphi} (x_0,0)\\-M\vert s\vert \partial_{x_0}\varphi(x_0,0) v_{\nu,\varphi}(x_0,0) -M\vert s\vert [\chi_\nu,\partial_{x_0}\varphi(x_0,0)]v_\varphi +r_{\nu,\varphi}.\nonumber
\end{eqnarray}
By (\ref{AAgora2A}) and (\ref{zara})
\begin{eqnarray}\label{zaraPP}
\int_0^T (\vert\partial^-_{x_1}v_{\nu,\varphi}\vert^2+\vert \tilde s\vert^3\vert v_{\nu,\varphi}\vert^2)(x_0,0)dx_0\le C_{45}(\vert \tilde s\vert\Vert P_{\varphi}(x,D,\widetilde s)v_{\nu,\varphi}\Vert^2_{L^2(Q_+)}\\
+\vert\tilde s\vert^2\varphi^{\frac {5}{6}}(x^*)\Vert v_\varphi(\cdot,0)\Vert^2_{L^2(0,T)}+\vert\tilde s\vert^2\Vert \partial^+_{x_1}v_{\varphi}(\cdot,0)\Vert^2_{L^2(0,T)}+ \Vert  r_{\nu,\varphi}\Vert^2_{L^2(0,T)}).\nonumber
\end{eqnarray}

Applying Proposition \ref{opana} to function $ v_{\nu,\varphi}$ on $Q_-$  and using (\ref{zaraPP}) to estimate the boundary integrals in (\ref{klop1}) we have

\begin{eqnarray}\label{zara2}
\vert\tilde s\vert\Vert v_{\nu,\varphi} \Vert^2_{H^{0,1,\tilde s}(Q_-)}\le C_{46}(\Vert P_{\varphi}(x,D,\widetilde s)v_{\nu,\varphi}\Vert^2_{Y}+\vert\tilde s\vert^3\varphi^{\frac {5}{6}}(x^*)\Vert v_\varphi(\cdot,0)\Vert^2_{L^2(0,T)}\nonumber\\
+\vert\tilde s\vert^2\Vert \partial^+_{x_1}v_{\varphi}(\cdot,0)\Vert^2_{L^2(0,T)}
+\vert\tilde s\vert\Vert  r_{\nu,\varphi}\Vert^2_{L^2(0,T)}).
\end{eqnarray}
Therefore (\ref{zara}) and  (\ref{zara2})  imply
\begin{eqnarray}\label{Pzara2L}
\vert\tilde s\vert\Vert v_{\nu,\varphi} \Vert^2_{H^{0,1,\tilde s}(Q_-)}+\vert\tilde s\vert^3\Vert v_{\nu,\varphi} \Vert^2_{H^{0,1,\tilde s}(Q_+)}+\Vert \mathcal  B v_{\nu,\varphi}\Vert^2_{\mathcal Z(0,T)}\le C_{47}(\Vert P_{\varphi}(x,D,\widetilde s)v_{\nu,\varphi}\Vert^2_{Y}\nonumber\\
+\vert\tilde s\vert^3\varphi^{\frac {5}{6}}(x^*)\Vert v_\varphi(\cdot,0)\Vert^2_{L^2(0,T)}+\vert\tilde s\vert\Vert \partial^+_{x_1}v_{\varphi}(\cdot,0)\Vert^2_{L^2(0,T)}
+\vert\tilde s\vert
\Vert r_{\nu,\varphi}\Vert^2_{L^2(0,T)}).
\end{eqnarray}

By (\ref{Pzara2L}) and (\ref{ladymikki}) we have
\begin{eqnarray}\label{zara2L}
\vert\tilde s\vert\Vert v_{\nu,\varphi} \Vert^2_{H^{\frac 12,1,\tilde s}(Q_-)}+\vert\tilde s\vert^3\Vert v_{\nu,\varphi} \Vert^2_{H^{\frac 12,1,\tilde s}(Q_+)}+\Vert \mathcal  B v_{\nu,\varphi}\Vert^2_{\mathcal Z(0,T)}\le C_{48}(\Vert P_{\varphi}(x,D,\widetilde s)v_{\nu,\varphi}\Vert^2_{Y}\nonumber\\
+\vert\tilde s\vert^3\varphi^{\frac {5}{6}}(x^*)\Vert v_\varphi(\cdot,0)\Vert^2_{L^2(0,T)}+\vert\tilde s\vert^2\Vert \partial^+_{x_1}v_{\varphi}(\cdot,0)\Vert^2_{L^2(0,T)}
+\vert\tilde s\vert
\Vert r_{\nu,\varphi}\Vert^2_{L^2(0,T)}).
\end{eqnarray}
{\bf Case 2.} Let $(x^*,\xi_0^*,\tilde s^*)\in \mathcal Z_{\varphi,-}(\ell)$ and $\mbox{supp}\, \chi_\nu\subset \mathcal O(\zeta^*,\delta_1).$ Then
$
\xi_0^*\ne 0
$ and decomposition (\ref{min}) holds true.
By Proposition \ref{gorokx} and (\ref{koma}) we have
\begin{eqnarray}\label{Gzopa1}
\Vert V^+_{\nu,\varphi}(\cdot,0)\Vert_{H^{\frac{1}{4}}(0,T)\cap L^2_{\tilde s}(0,T)}+\Vert V^+_{\nu,\varphi}\Vert_{H^{\frac{1}{2},1,\tilde s}(Q_+)}\\\le C_{49}(\Vert P_{\varphi}(x,D,\widetilde s)v_{\nu,\varphi}\Vert_{L^2(Q_+)}+\widetilde\varphi^\frac {5}{12}(x^*)\Vert v_\varphi\Vert_{H^{\frac 12,1,\tilde s}(Q_+)}).\nonumber
\end{eqnarray}
By (\ref{3.35}) for any $x\in B(x,\delta)\cap \mbox{supp}\, \eta_\ell$   and $(\xi,\tilde s)\in \mbox{supp}\,\chi_\nu$ we have  $\mbox{Im}\,\Gamma^+_\varphi(x,\xi,\tilde s)=-\vert\tilde s\vert\phi_0(x,x^*)+\frac{1}{\root\of{2}}\root\of {\rho(x
)\vert\xi_0\vert}.$
Since  $(x^*,\xi_0^*,\tilde s^*)\in \mathcal Z_{\varphi,-}(\ell)$ there exist  a positive $\delta_3=\delta_3(x^*,\zeta^*)$ such that
\begin{equation}\label{robinQ}
\vert \tilde s^*\vert \phi_0(x^*,x^*)- \frac{1}{\root\of{2}}\root\of {\rho(x^*
)\vert\xi_0^*\vert} <-\delta_3<0.
\end{equation}
Therefore there exists positive  $\delta_4(\delta_3) $ such that
\begin{equation}\label{robin}
\mbox{Im}\,\Gamma^\pm_\varphi(x, \xi,\tilde s ) >\delta_3/2 >0\quad (x,\xi,\tilde s)\in B(x,\delta_4)\cap \mbox{supp}\, \eta_\ell\times  \mbox{supp}\,\kappa(\nu,\cdot).
\end{equation}

By Proposition \ref{gorokxY} and (\ref{robin}) we have
\begin{equation}\label{Gzopa2}
\Vert U^-_{\nu,\varphi}(\cdot,0)\Vert_{H^{\frac{1}{4}}(0,T)\cap L^2_{ \tilde s}(0,T)}\le C_{50}(\Vert P_{\varphi}(x,D,\widetilde s)v_{\nu,\varphi}\Vert_{L^2(Q_-)}+\widetilde\varphi^\frac {5}{12}(x^*)\Vert v_\varphi\Vert_{H^{\frac 12,1,\tilde s}(Q_-)}).
\end{equation}

The equality (\ref{AAgora2A}) we replace $\partial^+_{x_1}v_{\nu,\varphi}(\cdot,0)$ by $i(\Gamma_\varphi^- (x,D_0,\tilde s)v_{\nu,\varphi}(\cdot,0)+V^+_{\nu,\varphi}(\cdot,0))$ and we replace  $\partial^-_{x_1}v_{\nu,\varphi}(\cdot,0)$ by $i(\Gamma_\varphi^+ (x,D_0,\tilde s) v_{\nu,\varphi}(\cdot,0)+U^-_{\nu,\varphi}(\cdot,0)).$

By (\ref{Gzopa1}) and (\ref{Gzopa2}) there exist a function $q$  such that on $[0,T]$ we have
\begin{eqnarray}\label{zopali}
M\partial_{x_0}v_{\nu,\varphi}(x_0,0)=i(\Gamma_\varphi^-(x,D_0,\tilde s) v_{\nu,\varphi}-\mu\Gamma_\varphi^+(x,D_0,\tilde s) v_{\nu,\varphi})(x_0,0)\\+M\vert s \vert \partial_{x_0}\varphi(x_0,0) v_{\nu,\varphi}(x_0,0) +q ,\nonumber
\end{eqnarray}
where
\begin{eqnarray}\label{zopa10}
\root\of{\vert \tilde s\vert}\Vert q\Vert_{L^2(0,T)}\le  C_{51}(\Vert P_{\varphi}(x,D,\widetilde s)v_{\nu,\varphi}\Vert_{L^2(Q)}+\root\of{\vert \tilde s\vert}\Vert  r_{\nu,\varphi}\Vert_{L^2(0,T)}\nonumber\\+\vert\tilde s\vert^\frac 43\Vert v_\varphi(\cdot,0)\Vert_{L^2(0,T)}+\root\of{\vert\tilde s\vert}\Vert \partial^+_{x_1}v_{\varphi}(\cdot,0)\Vert_{L^2(0,T)}+\widetilde\varphi^\frac {5}{12}(x^*)\Vert v_\varphi\Vert_{H^{\frac 12,1,\tilde s}(Q)}).
\end{eqnarray}

Inequalities  (\ref{robinQ}), (\ref{nikolas}) and short computations imply
\begin{eqnarray}\label{bigzopaP}
\Vert \partial_{x_0}  v_{\nu,\varphi}(\cdot,0)\Vert_{L^2(0,T)}\ge \Vert\eta_\ell \partial_{x_0} \chi_\nu v_{\varphi}(\cdot,0)\Vert_{L^2(0,T)}-\Vert \partial_{x_0}\eta_\ell  \chi_\nu v_{\varphi}(\cdot,0)\Vert_{L^2(0,T)}\ge\nonumber\\
\Vert \partial_{x_0} \chi_\nu v_{\varphi}(\cdot,0)\Vert_{L^2(0,T)}
-\Vert(1-\eta_\ell) \partial_{x_0} \chi_\nu v_{\varphi}(\cdot,0)\Vert_{L^2(0,T)}-\Vert \partial_{x_0}\eta_\ell  \chi_\nu v_{\varphi}(\cdot,0)\Vert_{L^2(0,T)}\ge\nonumber\\
 C_{52}\vert \tilde s\vert^2  \Vert  \chi_\nu v_{\varphi}(\cdot,0)\Vert_{L^2(\Bbb R^1)}-C_{53}\vert\tilde s\vert\varphi^\frac {5}{12}(x^*)\Vert v_\varphi\Vert_{L^2(0,T)}.
\end{eqnarray}
By (\ref{bigzopaP}), (\ref{norm})  and (\ref{zopa10}) we have
\begin{eqnarray}\label{xui}
\root\of{\vert \tilde s\vert}\Vert v_{\nu,\varphi}(\cdot,0)\Vert_{H^{1,\tilde s}(0,T)} +\vert\tilde s\vert^\frac 52 \Vert v_{\nu,\varphi}(\cdot,0)\Vert_{L^2(0,T)} \le  C_{54}(\Vert P_{\varphi}(x,D,\widetilde s)v_{\nu,\varphi}\Vert_{L^2(Q)}\nonumber\\+\root\of{\vert \tilde s\vert}\Vert r_{\nu,\varphi}\Vert_{L^2(0,T)}+\vert\tilde s\vert^2\widetilde\varphi^\frac {5}{12}(x^*)\Vert v_\varphi(\cdot,0)\Vert_{L^2(0,T)}+\widetilde\varphi^\frac {5}{12}(x^*)\Vert v_\varphi\Vert_{H^{\frac 12,1,\tilde s}(Q)}).
\end{eqnarray}

So  by (\ref{xui}) for any $\delta_3 >0$
\begin{eqnarray}\label{grom}
\vert \tilde s\vert^2\vert \mbox{Re}\,
\int_{\Bbb R^1} \partial_{x_1}^+ v_{\nu,\varphi}
\overline{\rho\partial_{x_0}v_{\nu,\varphi}}\vert_{x_1=0} d x_0\vert\quad\quad\quad\nonumber\\\le \delta_3\int_0^T (\vert \tilde s\vert\phi_0(x^*,x^*)\vert\partial^+_{x_1}v_{\nu,\varphi}\vert^2+\vert\tilde s\vert^3\phi_0^3(x^*,x^*)\vert v_{\nu,\varphi}\vert^2)(x_0,0)dx_0\nonumber\\+C_{55}(\delta_3)\vert\tilde s\vert\Vert \partial_{x_0}v_{\nu,\varphi}(\cdot,0)\Vert^2_{L^2(0,T)}\nonumber\\\le \delta_3\int_0^T (\vert \tilde s\vert\phi_0(x^*,x^*)\vert\partial^+_{x_1}v_{\nu,\varphi}\vert^2+\vert\tilde s\vert^3\phi_0^3(x^*,x^*)\vert v_{\nu,\varphi}\vert^2)(x_0,0)dx_0\nonumber\\
 +C_{56}(\delta_3)(\Vert P_{\varphi}(x,D,\widetilde s)v_{\nu,\varphi}\Vert^2_{L^2(Q)}+\vert \tilde s\vert\Vert r_{\nu,\varphi}\Vert^2_{L^2(0,T)}+\vert\tilde s\vert^2\widetilde\varphi^\frac {5}{6}(x^*)\Vert v_\varphi(\cdot,0)\Vert^2_{L^2(0,T)}\nonumber\\ +\widetilde\varphi^\frac {5}{6}(x^*)\Vert v_\varphi\Vert^2_{H^{\frac 12,1,\tilde s}(Q)}
 +\vert\tilde s\vert^2\Vert \partial^+_{x_1}v_{\varphi}(\cdot,0)\Vert^2_{L^2(0,T)}).
\end{eqnarray}
By Proposition \ref{opana} and estimate (\ref{grom})
\begin{eqnarray}\label{zaraP}
\int_0^T (\vert \tilde s\vert^3\vert\partial^+_{x_1}v_{\nu,\varphi}\vert^2+\vert\tilde  s\vert^5\vert v_{\nu,\varphi}\vert^2)(x_0,0)dx_0+\vert\tilde s\vert^3\Vert v_{\nu,\varphi} \Vert^2_{H^{0,1,\tilde s}(Q_+)}\le C_{57}(\Vert P_{\varphi}(x,D,\widetilde s)v_{\nu,\varphi}\Vert^2_{Y}\nonumber\\+\vert\tilde s\vert\Vert r_{\nu,\varphi}\Vert^2_{L^2(0,T)}+\vert\tilde s\vert^2\Vert \partial^+_{x_1}v_{\varphi}(\cdot,0)\Vert^2_{L^2(0,T)}+\vert\tilde s\vert^2\widetilde\varphi^\frac {5}{6}(x^*)\Vert v_\varphi(\cdot,0)\Vert^2_{L^2(0,T)}\nonumber\\+\widetilde\varphi^\frac {5}{6}(x^*)\Vert v_\varphi\Vert^2_{H^{\frac 12,1,\tilde s}(Q)}).
\end{eqnarray}
Hence
\begin{eqnarray}\label{zopa4}
\int_0^T (\vert\tilde s\vert \vert\partial^-_{x_1}v_{\nu,\varphi}\vert^2+\vert\tilde  s\vert^5 \vert v_{\nu,\varphi}\vert^2)(x_0,0)dx_0
\le C_{58}(\Vert P_{\varphi}(x,D,\widetilde s)v_{\nu,\varphi}\Vert^2_{Y}+\vert\tilde s\vert\Vert  r_{\nu,\varphi}\Vert^2_{L^2(0,T)}\nonumber\\+\vert\tilde s\vert^2\widetilde\varphi^\frac {5}{6}(x^*)\Vert v_\varphi(\cdot,0)\Vert^2_{L^2(0,T)}+\vert\tilde s\vert^2\Vert \partial^+_{x_1}v_{\varphi}(\cdot,0)\Vert^2_{L^2(0,T)}+\widetilde\varphi^\frac {5}{6}(x^*)\Vert v_\varphi\Vert^2_{H^{\frac 12,1,\tilde s}(Q)}).
\end{eqnarray}
By (\ref{zopa4}) and inequality (\ref{klop1}) of Proposition \ref{opana}
\begin{eqnarray}\label{zara1L}
\vert\tilde s\vert\Vert v_{\nu,\varphi} \Vert^2_{H^{0,1,\tilde s}(Q_-)}\le C_{59}(\Vert P_{\varphi}(x,D,\widetilde s)v_{\nu,\varphi}\Vert^2_{Y}+\vert \tilde s\vert\Vert r_{\nu,\varphi}\Vert^2_{L^2(0,T)}\nonumber\\+\vert\tilde s\vert^2\widetilde\varphi^\frac {5}{6}(x^*)\Vert v(\cdot,0)\Vert^2_{L^2(0,T)}+\vert\tilde s\vert^2\Vert \partial^+_{x_1}v_{\varphi}(\cdot,0)\Vert^2_{L^2(0,T)}+\widetilde\varphi^\frac {5}{6}(x^*)\Vert v_\varphi\Vert^2_{H^{\frac 12,1,\tilde s}(Q)}).
\end{eqnarray}
 By (\ref{Gzopa1}), (\ref{min11}) and  (\ref{zaraP}) the following estimate is true
 \begin{eqnarray}\label{fool}
 \vert \tilde s\vert^3\Vert \alpha^+(x,D_0,\tilde s)v_{\nu,\varphi}\Vert^2_{L^2(Q_+)}+\vert \tilde s\vert\Vert \alpha^+(x,D_0,\tilde s)v_{\nu,\varphi}\Vert^2_{L^2(Q_-)}\\\le  C_{60}(\Vert P_{\varphi}(x,D,\widetilde s)v_{\nu,\varphi}\Vert^2_{Y}+\vert\tilde s\vert\Vert r_{\nu,\varphi}\Vert^2_{L^2(0,T)}+\vert\tilde s\vert^2\Vert \partial^+_{x_1}v_{\varphi}(\cdot,0)\Vert^2_{L^2(0,T)}\nonumber\\+\Vert s\vert^2\widetilde\varphi^\frac {5}{6}(x^*)\Vert v_\varphi(\cdot,0)\Vert^2_{L^2(0,T)}+\widetilde\varphi^\frac {5}{6}(x^*)\Vert v_\varphi\Vert^2_{H^{\frac 12,1,\tilde s}(Q)}).\nonumber
 \end{eqnarray}
 Using the G\aa rding inequality proved in  we obtain from (\ref{fool})
 \begin{eqnarray}\label{fool1}
 \vert \tilde s\vert^3\Vert \alpha^+(x,D_0,\tilde s)v_{\nu,\varphi}\Vert^2_{H^{\frac 12,0}(Q_+)}\le  C_{61}(\Vert P_{\varphi}(x,D,\widetilde s)v_{\nu,\varphi}\Vert^2_{Y}+\vert\tilde s\vert\Vert  r_{\nu,\varphi}\Vert^2_{L^2[0,T]}\\+\vert\tilde s\vert^2\Vert \partial^+_{x_1}v_{\varphi}(\cdot,0)\Vert^2_{L^2(0,T)}+\vert\tilde s\vert^2\widetilde\varphi^\frac {5}{6}(x^*)\Vert v_\varphi(\cdot,0)\Vert^2_{L^2(0,T)}\nonumber+\widetilde\varphi^\frac {5}{6}(x^*)\Vert v_\varphi\Vert^2_{H^{\frac 12,1,\tilde s}(Q)}).
 \end{eqnarray}
The inequality  (\ref{zara1L}), (\ref{fool1}) and  (\ref{zaraP}) we have

\begin{eqnarray}\label{zara2LL}
\vert\tilde s\vert\Vert v_{\nu,\varphi} \Vert^2_{H^{\frac 12,1,\tilde s}(Q_-)}+\vert\tilde s\vert^3\Vert v_{\nu,\varphi} \Vert^2_{H^{\frac 12,1,\tilde s}(Q_+)}+\Vert \mathcal  B v_{\nu,\varphi}\Vert^2_{\mathcal Z(0,T)}\le C_{62}(\Vert P_{\varphi}(x,D,\widetilde s)v_{\nu,\varphi}\Vert^2_{Y}\\
+\vert\tilde s\vert^3\varphi^{\frac {5}{6}}(x^*)\Vert v_\varphi(\cdot,0)\Vert^2_{L^2(0,T)}+\vert\tilde s\vert^2\Vert \partial^+_{x_1}v_{\varphi}(\cdot,0)\Vert^2_{L^2(0,T)}
+\vert\tilde s\vert
\Vert r_{\nu,\varphi}\Vert^2_{L^2(0,T)}+\widetilde\varphi^\frac {5}{6}(x^*)\Vert v_\varphi\Vert^2_{H^{\frac 12,1,\tilde s}(Q)})\nonumber.
\end{eqnarray}
{\bf Case 3.} {\it  Let $(x^*,\xi_0^*,\tilde s^*)\in \mathcal Z_{\varphi,0}(\ell)$ and $ \mbox{supp}\, \chi_\nu\subset \mathcal O(\zeta^*,\delta_1).$} Consider two subcases:

{\bf Subcase 1.} {\it Let $
(x^*,\xi_0^*,\tilde s^*)\in \mathcal Z_{\varphi_1,0}(\ell).$}  Then by (\ref{gopnikz}) and (\ref{barbos}) $
(x^*,\xi_0^*,\tilde s^*)\in \mathcal Z_{\varphi_2,+}(\ell).$  Therefore estimate (\ref{zara}) is true:
\begin{eqnarray}\label{Zzara}
\int_0^T (\vert \tilde s\vert \vert\partial^+_{x_1}v_{\nu ,\varphi_2}\vert^2+\vert \tilde s\vert^3\vert v_{\nu,\varphi_2}\vert^2)(x_0,0)dx_0+\vert\tilde s\vert\Vert v_{\nu ,\varphi_2} \Vert^2_{H^{0,1,\tilde s}(Q_+)}\\\le C_{63}(\Vert P_{\varphi_2}(x,D,\widetilde s)v_{\nu,\varphi_2} \Vert^2_{L^2(Q_+)}+\vert\tilde s\vert\tilde\varphi^\frac 56(x^*)\Vert v_\varphi(\cdot,0)\Vert^2_{L^2(0,T)}).\nonumber
\end{eqnarray}
Since by (\ref{barbos}) on  the segment  $[0,T]$
$
\varphi_1(x_0,0)=\varphi_2(x_0,0)
$ the inequality (\ref{Zzara}) implies
\begin{eqnarray}\label{Zzara0}
\int_0^T (\vert\tilde  s\vert^3\vert\partial^+_{x_1}v_{\nu ,\varphi_1}\vert^2+\vert \tilde s\vert^5\vert v_{\nu ,\varphi_1}\vert^2)(x_0,0)dx_0+\vert\tilde s\vert^3\Vert v_{\nu,\varphi_2} \Vert^2_{H^{0,1,\tilde s}(Q_+)}\\\le C_{64} (\vert\tilde  s\vert^2\Vert P_{\varphi_2}(x,D,\widetilde s)v_{\nu,\varphi_2} \Vert^2_{L^2(Q_+)}+\vert\tilde s\vert^3\tilde\varphi^\frac 56(x^*)\Vert v_\varphi(\cdot,0)\Vert^2_{L^2(0,T)}).\nonumber
\end{eqnarray}
Observe that we have  equality (\ref{AAgora2A}) with function $\varphi_1$. Using this equality and (\ref{Zzara00}) we  estimate $\partial^-_{x_1}v_{\nu ,\varphi_1}(\cdot,0)$:

\begin{eqnarray}\label{LAAgora2A}
\Vert\partial^-_{x_1}v_{\nu,\varphi_1} (\cdot,0)\Vert_{L^2(0,T)}\le \Vert\partial^+_{x_1}v_{\nu,\varphi_1} (\cdot,0)\Vert_{L^2(0,T)}+M\Vert\partial_{x_0}v_{\nu,\varphi_*} (\cdot,0)\Vert_{L^2(0,T)}\\+M\Vert\vert s\vert \partial_{x_0}\varphi_* v_{\nu,\varphi_*}(\cdot,0) \Vert_{L^2(0,T)} +M\vert s\vert \Vert[\chi_\nu,\partial_{x_0}\varphi_*(x_0,0)]v_{\varphi_*}\Vert_{L^2(0,T)} +\Vert \tilde r_{\nu,\varphi_*}\Vert_{L^2(0,T)}\nonumber\\
+C_{65}(\vert\tilde s\vert\tilde\varphi^{\frac {5}{12}}(x^*)\Vert v_{\varphi_*}(\cdot,0)\Vert_{L^2(0,T)}+
\vert \tilde s\vert^2\Vert v_{\nu,\varphi_*}(\cdot,0)\Vert_{L^2(\Bbb R^1)}) .\nonumber
\end{eqnarray}

Applying the Proposition \ref{opana}  and  using (\ref{LAAgora2A}) and (\ref{Zzara00}) to estimate the boundary integrals in  the right hand side of (\ref{klop1}) we have
\begin{equation}\label{Zzara2z}
 \vert\tilde s\vert\Vert v_{\nu,\varphi_1} \Vert^2_{H^{0,1,\tilde s}(Q_-)}\le C_{66}(\Vert P_{\varphi_*}(x,D,\widetilde s)v_{\nu,\varphi_*}\Vert^2_{Y}
+\vert\tilde s\vert^3\tilde\varphi^\frac 56(x^*)\Vert v_{\varphi_*}(\cdot,0)\Vert^2_{L^2(0,T)}
+\vert \tilde s\vert \Vert r_{\nu,\varphi_*}\Vert^2_{L^2(0,T)}).
\end{equation}
Combining the estimates (\ref{zaraPP}) (\ref{Zzara00}) and (\ref{Zzara2z}) we have (\ref{zara2LL}).
Now  we consider

{\bf Subcase 2.} {\it Let   $
(x^*,\xi_0^*, \tilde s^*)\in \mathcal Z_{\varphi_2,0}(\ell)\quad \mbox{and}\quad \mbox{supp}\, \chi_\nu\subset \mathcal O(\zeta^*,\delta_1).$}  Then $
(x^*,\xi_0^*,\tilde s^*)\in \mathcal Z_{\varphi_1,-}(\ell)$ and  inequality (\ref{robinQ}) holds true for some positive $\delta_3.$
Therefore there exists positive  $\delta_4(\delta_3) $ such that
\begin{equation}\label{robin3}
\mbox{Im}\,\Gamma^\pm_{\varphi_1}(x,\xi,\tilde s ) >\delta_3/2>0 \quad (x,\xi,\tilde s)\in B(x^*,\delta_4)\cap \mbox{supp}\, \eta_\ell\times  \mbox{supp}\,\kappa(\nu,\cdot).
\end{equation}

By Proposition \ref{gorokx} and (\ref{robin3}) we have
\begin{equation}\label{Lzopa2}
\Vert U^-_{\nu ,\varphi_1}(\cdot,0)\Vert_{H^{\frac{1}{4}}(0,T)\cap L^2_{ \tilde s}(0,T)}\le C_{67}(\Vert P_{\varphi_1}(x,D,\widetilde s)v_{\nu,\varphi_1}\Vert_{L^2(Q_-)}+\widetilde\varphi^\frac {5}{12}(x^*)\Vert v_{\varphi_1}\Vert_{H^{\frac 12,1,\tilde s}(Q_-)}).
\end{equation}
By Proposition \ref{gorokx} we have
\begin{equation}\label{Lzopa1}
\Vert V^+_{\nu,\varphi_2}(\cdot,0)\Vert_{H^{\frac{1}{4}}(0,T)\cap L^2_{\tilde s}(0,T)}\le C_{68}(\Vert P_{\varphi_2}(x,D,\widetilde s)v_{\nu,\varphi_2}\Vert_{L^2(Q_+)}+\widetilde\varphi^\frac {5}{12}(x^*)\Vert v_{\varphi_2}\Vert_{H^{\frac 12,1,\tilde s}(Q_+)}).
\end{equation}

The equality (\ref{gora2A}) for the function $u$ can be written as
 \begin{eqnarray}\label{Agora2AL}
-\partial^-_{x_1}v_{\varphi_1}(x_0,0)+\mu\partial^+_{x_1}v_{\varphi_2}(x_0,0)+\vert s\vert(\partial_{x_1}\varphi_1 -\mu\partial_{x_1}\varphi_2) v_{\varphi_*}=\\ \nonumber M\partial_{x_0}v_{\varphi_*}(x_0,0)-M\vert s\vert \partial_{x_0}\varphi_*(x_0,0) v_{\varphi_*}(x_0,0)+ r_{\varphi_*}.
\end{eqnarray}

We apply to both sides of equation (\ref{Agora2AL}) the operator $\chi_\nu(x,D_0,\tilde s):$
\begin{eqnarray}\label{AAgora2AL}
-\partial^-_{x_1}v_{\nu,\varphi_1}(x_0,0)+\mu\partial^+_{x_1}v_{\nu,\varphi_2}(x_0,0)+\vert s\vert(\partial_{x_1}\varphi_1 -\mu\partial_{x_1}\varphi_2) v_{\nu,\varphi_*}(x_0,0)\\ \nonumber-M\partial_{x_0}v_{\nu,\varphi_*}(x_0,0)+M\vert s\vert \partial_{x_0}\varphi_*(x_0,0) v_{\nu,\varphi_*}(x_0,0)- r_{\nu,\varphi_*}(x_0)\nonumber\\
-[\chi_\nu,\partial^-_{x_1}]v_{\varphi_1}(x_0,0)+[\chi_\nu,\mu\partial^+_{x_1}]v_{\varphi_2}(x_0,0)+\vert s\vert[\chi_\nu,(\partial_{x_1}\varphi_1 -\mu\partial_{x_1}\varphi_2)] v_{\varphi_*}(x_0,0)\nonumber\\ +M\nonumber[\chi_\nu,\partial_{x_0}]v_{\varphi_*}(x_0,0)+M\vert s\vert [\chi_\nu,\partial_{x_0}\varphi_*(x_0,0)] v_{\varphi_*}(x_0,0).
\end{eqnarray}

In the equality (\ref{AAgora2AL}) we replace $\partial^+_{x_1}v_{\nu,\varphi_2}(\cdot,0)$ by $i(\Gamma_{\varphi_2}^- (x,D_0,\tilde s) v_{\nu,\varphi_2}(\cdot,0)+V^+_{\nu,\varphi_2}(\cdot,0))$ and we replace  $\partial^-_{x_1}v_{\nu,\varphi_1}(\cdot,0)$ by $i(\Gamma_{\varphi_1}^+(x,D_0,\tilde s) v_{\nu,\varphi_1}(\cdot,0)+U^-_{\nu,\varphi_1}(\cdot,0)):$
\begin{eqnarray}\label{zombie}
M\partial_{x_0}v_{\nu ,\varphi_*}(x_0,0)=i(\mu\Gamma_{\varphi_2}^-(x,D_0,\tilde s) v_{\nu ,\varphi_2}-\Gamma_{\varphi_1}^+ (x,D_0,\tilde s)v_{\nu ,\varphi_1})(x_0,0)\nonumber\\+i(\mu V^+_{\nu,\varphi_2}(x_0,0)- U^-_{\nu,\varphi_1}(x_0,0))\nonumber\\
-\vert s\vert(\partial_{x_1}\varphi_1 -\mu\partial_{x_1}\varphi_2) v_{\nu,\varphi_*}(x_0,0)-M\vert s\vert \partial_{x_0}\varphi_*(x_0,0) v_{\nu,\varphi_*}(x_0,0)+ r_{\nu,\varphi_*}(x_0)\nonumber\\
+[\chi_\nu,\partial^-_{x_1}]v_{\varphi_1}(x_0,0)-[\chi_\nu,\mu\partial^+_{x_1}]v_{\varphi_2}(x_0,0)-\vert s\vert[\chi_\nu,(\partial_{x_1}\varphi_1 +\mu\partial_{x_1}\varphi_2)] v_{\varphi_*}(x_0,0)\nonumber\\ -M[\chi_\nu,\partial_{x_0}]v_{\varphi_*}(x_0,0)-M\vert s\vert [\chi_\nu,\partial_{x_0}\varphi_*(x_0,0)] v_{\varphi_*}(x_0,0).
\end{eqnarray}

Observe that estimate (\ref{bigzopaP}) holds true.
By (\ref{bigzopaP}), (\ref{Lzopa2}), (\ref{Lzopa1}) from (\ref{zombie}) we have
\begin{eqnarray}\label{Axuii}
\root\of{\vert\tilde s\vert}\Vert v_{\nu ,\varphi_*}(\cdot,0)\Vert_{H^{1,\tilde s}(0,T)}+\vert\tilde s\vert^\frac 52\Vert v_{\nu ,\varphi_*}(\cdot,0)\Vert_{L^2(0,T)}\le   C_{69}(\root\of{\vert\tilde s\vert}\Vert P_{\varphi_*}(x,D,\widetilde s)v_{\nu,\varphi_*}\Vert_{L^2(Q)}\nonumber\\+\widetilde\varphi^\frac {5}{12}(x^*)\Vert v_{\varphi_*}\Vert_{H^{\frac 12,1,\tilde s}(Q)}+\root\of{\vert\tilde s\vert}\Vert  r_{\nu,\varphi_*}\Vert_{L^2(0,T)}+\vert\tilde s\vert^\frac 32\tilde\varphi^\frac{5}{12}(x^*)\Vert v_{\varphi_*}(\cdot,0)\Vert_{L^2(0,T)}).
\end{eqnarray}.
By (\ref{Axuii}) and  Proposition \ref{opana} we obtain
\begin{eqnarray}\label{dondon}
\int_0^T (\vert\tilde s\vert^3 \vert\partial^+_{x_1}v_{\nu ,\varphi_2}\vert^2+\vert \tilde s\vert^5 \vert v_{\nu ,\varphi_2}\vert^2)(x_0,0)dx_0+\vert\tilde s\vert^3\Vert v_{\nu,\varphi_2} \Vert^2_{H^{0,1,\tilde s}(Q_+)}\\\le C_{70}(\Vert P_{\varphi_*}(x,D,\widetilde s)v_{\nu,\varphi_*}\Vert^2_{Y}+\widetilde\varphi^\frac {5}{6}(x^*)\Vert v_{\varphi_*}\Vert^2_{H^{\frac 12,1,\tilde s}(Q)}+\vert \tilde s\vert\Vert \tilde r_{\nu,\varphi_*}\Vert^2_{L^2(0,T)}\nonumber\\+\vert\tilde s\vert^2\Vert \partial^+_{x_1}v_{\varphi_*}(\cdot,0)\Vert^2_{L^2(0,T)}+\vert\tilde s\vert^3 \tilde\varphi^\frac{5}{6}(x^*)\Vert v_{\varphi_*}(\cdot,0)\Vert^2_{L^2(0,T)}).\nonumber
\end{eqnarray}
Hence
\begin{eqnarray}\label{Pzopa4}
\int_0^T (\vert\tilde s\vert \vert\partial^-_{x_1}v_{\nu,\varphi_2}\vert^2+\vert\tilde s\vert^3\vert v_{\nu,\varphi_2}\vert^2)(x_0,0)dx_0
\le C_{71}(\Vert P_{\varphi_*}(x,D,\widetilde s)v_{\nu,\varphi_*}\Vert^2_{L^2(Q)}\\+\widetilde\varphi^\frac {5}{6}(x^*)\Vert v_{\varphi_*}\Vert^2_{H^{\frac 12,1,\tilde s}(Q)}+\vert\tilde s\vert\Vert r_{\nu,\varphi_*}\Vert^2_{L^2(0,T)}+\vert\tilde s\vert^2\Vert \partial^+_{x_1}v_{\varphi}(\cdot,0)\Vert^2_{L^2(0,T)}\nonumber\\+\vert\tilde s\vert^3\tilde \varphi^\frac 56(x^*)\Vert v_{\varphi_*}(\cdot,0)\Vert^2_{L^2(0,T)}).\nonumber
\end{eqnarray}
By Proposition \ref{opana} and  estimate (\ref{Pzopa4})
\begin{eqnarray}\label{Pzara1L}
\vert\tilde s\vert\Vert v_{\nu,\varphi_2} \Vert^2_{H^{0,1,\tilde s}(Q_-)}\le C_{72}(\Vert P_{\varphi_*}(x,D,\widetilde s)v_{\nu,\varphi_*}\Vert^2_{L^2(Q)}+\widetilde\varphi^\frac {5}{6}(x^*)\Vert v_{\varphi_*}\Vert^2_{H^{\frac 12,1,\tilde s}(Q)}\\+\vert \tilde s\vert\Vert r_{\nu,\varphi_*}\Vert^2_{L^2(0,T)}+\vert\tilde s\vert^2\Vert \partial^+_{x_1}v_{\varphi_*}(\cdot,0)\Vert^2_{L^2(0,T)}+\vert\tilde s\vert^3\tilde\varphi^\frac 56(x^*)\Vert v_{\varphi_*}(\cdot,0)\Vert^2_{L^2(0,T)}).\nonumber
\end{eqnarray}

The estimates (\ref{zaraPP}), (\ref{Pzara1L}) and (\ref{dondon}) imply
\begin{eqnarray}\label{zara2LLL}\int_0^T (\vert\tilde s\vert^3 \vert\partial^+_{x_1}v_{\nu ,\varphi_2}\vert^2+\vert \tilde s\vert^5 \vert v_{\nu ,\varphi_2}\vert^2)(x_0,0)dx_0+\vert\tilde s\vert^3\Vert v_{\nu,\varphi_2} \Vert^2_{H^{0,1,\tilde s}(Q_+)}\\
+\vert\tilde s\vert\Vert v_{\nu,\varphi_2} \Vert^2_{H^{0,1,\tilde s}(Q_-)}\le C_{73}(\Vert P_{\varphi_*}(x,D,\widetilde s)v_{\nu,\varphi_*}\Vert^2_{L^2(Q)}+\widetilde\varphi^\frac {5}{6}(x^*)\Vert v_{\varphi_*}\Vert^2_{H^{\frac 12,1,\tilde s}(Q)}\nonumber\\+\vert \tilde s\vert\Vert r_{\nu,\varphi_*}\Vert^2_{L^2(0,T)}+\vert\tilde s\vert^2\Vert \partial^+_{x_1}v_{\varphi_*}(\cdot,0)\Vert^2_{L^2(0,T)}+\vert\tilde s\vert^3\tilde\varphi^\frac 56(x^*)\Vert v_{\varphi_*}(\cdot,0)\Vert^2_{L^2(0,T)}).\nonumber
\end{eqnarray}

Now we observe that in all three cases the estimate (\ref{zara2LLL}) is true.
In order to finish the proof of the Proposition \ref{zoopa}
 let us take the covering of the surface $\Bbb M=\{ (\xi_0,\tilde s)\vert M(\xi_0,\tilde s)=1\}$ by
conical neighborhoods $\mathcal O(\zeta^*,\delta_1(\zeta^*)).$
From this covering we take the finite subcovering $\cup_{\nu=1}^N\mathcal O
(\zeta_\nu^*,\delta_1(\zeta_\nu^*))$.  Let us show that such a subcovering can be taken independently of parameter $\ell$ for all $\ell\ge \ell_0.$
Indeed
by (\ref{zima}) for any $x^*_\ell$  from $\mbox{supp}\, \kappa_\ell$  then \begin{equation}\label{molnia} supp_{y\in\mbox{supp}\, \kappa_\ell} \vert \phi_0(y,x^*_\ell)-1\vert \rightarrow 0\quad\mbox{as}\,\,\ell\rightarrow +\infty.
\end{equation}
We set
$$\mathcal Z_{\pm}=\{(x,\xi_0,\tilde s)\in \{(0,0), (T,0)\} \times \Bbb M\vert \quad \pm( \vert\tilde s\vert-\mbox{Im}\, \root \of {i\rho(x_0,0)\xi_0}) >0 ,\,\, x_0\in \{0,T\}\}
$$
and
$$\mathcal Z_{0}=\{(x,\xi_0,\tilde s)\in \{(0,0), (T,0)\} \times \Bbb M\vert\quad  \vert\tilde s\vert=\mbox{Im}\, \root \of {i\rho(x_0,0)\xi_0} =0,\,\, x_0\in \{0,T\} \}.
$$
By (\ref{molnia})
\begin{equation}\label{zacada}
\mbox{dist}\, (\mathcal Z_{\pm},\mathcal Z_{\varphi,\pm}(\ell))+\mbox{dist}\, (\mathcal Z_{0},\mathcal Z_{\varphi,\pm}(\ell))\rightarrow 0\quad\mbox{as}\,\,\ell\rightarrow +\infty.
\end{equation}
Without loss of generality one may assume that $x^*=(0,0)$ or $x^*=(T,0)$ in (\ref{3.55}). Let $ x^*=(T,0).$   We construct the covering of the set $\Bbb M$ in the following way: if $ (\xi_0^*,\tilde s^*) \in \mathcal M_{\pm}=\{(\xi_0,\tilde s)\in  \Bbb M\vert \,\pm( \vert\tilde s\vert-\mbox{Im}\, \root \of {i\rho(T,0)\xi_0}) >0 \}$ we consider covering of this point  by the ball of centered at  $ (\xi_0^*,\tilde s^*)$ of sufficiently small  radius $\delta(\xi_0^*,\tilde s^*)$ such that $\overline{\Bbb M\cap B((\xi_0^*,\tilde s^*),\delta(\xi_0^*,\tilde s^*))}\subset \mathcal M_{\pm}.$  If $ (\xi_0^*,\tilde s^*) \in \mathcal M_{0}=\{(\xi_0,\tilde s)\in  \Bbb M\vert \,  \vert\tilde s\vert-\mbox{Im}\, \root \of {i\rho(T,0)\xi_0} =0 \}$ we consider covering of this point  by the ball of centered at  $ (\xi_0^*,\tilde s^*)$ of sufficiently small  radius $\delta(\xi_0^*,\tilde s^*). $ From this covering we take the finite subcovering $\cup_{\nu=1}^N\mathcal O
(\zeta_\nu^*,\delta_1(\zeta_\nu^*))$ and let $\chi_\nu$ be the partition of unity subjected to this finite subcovering. Hence $\sum_{\nu=1}^N\chi_\nu
( \xi_0,\widetilde s)\equiv 1$ or all $( \xi_0,\widetilde s)$ such that $M(\xi_0,\widetilde s)\ge 1.$
By (\ref{zacada}) there exists $\ell_1$ such that for all $\ell\ge \ell_1$ if $(\xi_0^*,\tilde s^*) \in \mathcal M_{\pm}$ and $\supp \chi_\nu\in\Bbb M\cap B((\xi_0^*,\tilde s^*),\delta(\xi_0^*,\tilde s^*))$ then $$
\mbox{supp}\, \kappa_\ell\times \{0\} \times \mbox{supp}\, \chi_\nu\subset \mathcal Z_{\varphi,\pm}(\ell).
$$
On the other hand if  $(\xi_0^*,\tilde s^*) \in \mathcal M_{0}$ and $\supp \chi_\nu\in\Bbb M\cap B((\xi_0^*,\tilde s^*),\delta(\xi_0^*,\tilde s^*))$ then for all sufficiently large $\ell$ there exists $ \zeta^*(\ell)=(\xi_0^*(\ell),\tilde s^*(\ell)) $ such that $\mbox{supp}\, \chi_\nu\subset \mathcal O(\zeta^*(\ell),\delta_3) $ where $\delta_3\rightarrow +0$ as $\delta(\xi_0^*,\tilde s^*)\rightarrow +0.$

For index $\ell\in \{1,\dots, \ell_0\}$ individual subcovering of the set $\Bbb M$ and the corresponding partition of unity $\chi_\nu.$
   Let $\chi_0( \xi_0,\widetilde s)\in C^\infty_0(\Bbb R^{2})$ be a nonnegative function  which is identically equal one  if $M( \xi_0,\widetilde s)\le 1.$
 Then by (\ref{zara2L}) we have
\begin{eqnarray}\label{golifax}
\Vert \mathcal B v_{\varphi_*}\Vert_{\mathcal Z(0,T)}
+\root\of{ \vert\widetilde s\vert} \Vert
v_{\varphi_*}\Vert_{H^{\frac 12,1,\widetilde s}(Q)}
\le C_{74}\sum_{\nu=0}^N ( \Vert \eta_\ell\chi_\nu \mathcal B v_{\varphi_*} \Vert_{\mathcal Z(\Bbb R^1)}
+ \root\of{\vert \widetilde s\vert}\Vert\eta_\ell
\chi_\nu v_{\varphi_*}\Vert_{H^{\frac 12,1,\widetilde s}( Q)}\nonumber\\
 +\Vert (1-\eta_\ell)\chi_\nu \mathcal B v_{\varphi_*} \Vert_{\mathcal Z(\Bbb R^1)}
+ \root\of{\vert \widetilde s\vert}\Vert(1-\eta_\ell)
\chi_\nu v_{\varphi_*}\Vert_{H^{\frac 12,1,\widetilde s}( Q)}) \le C_{75}(\widetilde \varphi^\frac{5}{12}(x^*)\Vert v_{\varphi_*}
\Vert_{H^{\frac 12,1,\widetilde s}( Q)}\nonumber\\
+ \sum_{\nu=0}^N\Vert \mbox{\bf P}(x,D,\widetilde s)v_{\nu,\varphi_*}\Vert
_{Y\times L^2_{\tilde s}(0,T)}
+ \vert\tilde s\vert^\frac 32 \tilde \varphi^\frac{5}{12}(x^*)\Vert v_{\varphi_*}(\cdot,0)
)\Vert_{L^2 (0,T)}+\vert\tilde s\vert\Vert \partial_{x_1}^+v_{\varphi_*}(\cdot,0)
)\Vert_{L^2 (0,T)}                   \nonumber\\
+\Vert (1-\eta_\ell)\chi_\nu\mathcal B  v_{\varphi_*}\Vert_{\mathcal Z(\Bbb R^1)}
+ \root\of{\vert \widetilde s\vert}\Vert(1-\eta_\ell)
\chi_\nu v_{\varphi_*}\Vert_{H^{\frac 12,1,\widetilde s}(Q)}) .
\end{eqnarray}
By (\ref{po1}) and (\ref{nikolasZ}) there exist a constant $C_{76}$ independent of $\widetilde s,\ell$ and $\nu$ such that
\begin{eqnarray}\label{nikolas211}
\sum_{\nu=0}^N (\Vert (1-\eta_\ell)\chi_\nu\mathcal B  v_{\varphi_*}\Vert_{\mathcal Z(\Bbb R^1)}
+ \root\of{\vert \widetilde s\vert}\Vert(1-\eta_\ell)\chi_\nu
v_{\varphi_*}\Vert_{H^{\frac 12,1,\widetilde s}( Q)})\nonumber\\
\le
C_{76}\left(\vert\tilde s\vert^\frac 32 \tilde \varphi^\frac{5}{12}(x^*)\Vert v_{\varphi_*}(\cdot,0)
)\Vert_{L^2 (0,T)}
+ \root\of{\widetilde \varphi(x^*)}\Vert
v_{\varphi_*}\Vert_{H^{\frac 12,1,\widetilde s}( Q)}\right).
\end{eqnarray}
Using  inequality (\ref{nikolas211}) in order to estimate the last terms in
(\ref{golifax}) we obtain
\begin{eqnarray}
\Vert \mathcal B v_{\varphi_*}\Vert_{\mathcal Z(0,T)}
+\root\of{ \vert\widetilde s\vert} \Vert
v_{\varphi_*}\Vert_{H^{\frac 12,1,\widetilde s}(Q)}\le C_{77}\biggl(\root\of{\widetilde \varphi(x^*)}\Vert v_{\varphi_*}
\Vert_{H^{\frac 12,1,\widetilde s}(Q)}
+\Vert \mbox{\bf P}(x,D,\widetilde s)
v_{\varphi_*}\Vert_{Y\times L^2_{\tilde s}(0,T)}\nonumber\\
+ \vert\tilde s\vert^\frac 32 \tilde \varphi^\frac{5}{12}(x^*)\Vert v_{\varphi_*}(\cdot,0)
)\Vert_{L^2 (0,T)}
+ \sum_{\nu=0}^N\Vert[\tilde\chi_\nu,\mbox{\bf  P}(x,D,\widetilde s)]v_{\varphi_*}\Vert_{Y\times L^2_{\tilde s}(0,T)}
\biggr)         \nonumber\\
\le C_{78}\biggl(\root\of{\widetilde \varphi(x^*)}\Vert v_{\varphi_*}
\Vert_{H^{\frac 12,1,\widetilde s}( Q)}
+ \Vert\mbox{\bf  P}(x,D,\widetilde s)v_{\varphi_*}\Vert_{Y\times L^2_{s\tilde \varphi}(0,T)}
+ \vert\tilde s\vert^\frac 32 \tilde \varphi^\frac{5}{12}(x^*)\Vert v_{\varphi_*}(\cdot,0)
)\Vert_{L^2 (0,T)}\biggr)\nonumber.
\end{eqnarray}
Hence there exists $s_0>1$  such that for all
$s\ge s_0$  we see
\begin{eqnarray}\label{zmpolit}
\Vert \mathcal B v_{\varphi_*}\Vert_{\mathcal Z(0,T)}
+\root\of{ \vert\widetilde s\vert} \Vert
v_{\varphi_*}\Vert_{H^{\frac 12,1,\widetilde s}(Q)}
                         \\
\le C_{79}\biggl(\root\of{\widetilde \varphi(x^*)}\Vert v_{\varphi_*}
\Vert_{H^{\frac 12,1,\widetilde s}(Q)}
+ \Vert \mbox{\bf P}(x,D,\widetilde s)v_{\varphi_*}\Vert_{Y\times L^2_{s\tilde \varphi}(0,T)}+ \vert\tilde s\vert^\frac 32 \tilde \varphi^\frac{5}{12}(x^*)\Vert v_{\varphi_*}(\cdot,0)
)\Vert_{L^2 (0,T)} \biggr).\nonumber
\end{eqnarray}
Proof of Proposition \ref{zoopa} is complete. $\blacksquare$

 We set  \begin{equation}\label{lodka}\psi^*(x)=\left \{\begin{matrix}\hat s\varphi_*\circ F^{-1}(x)\quad\mbox{for}\,\,x_0\in [\frac T2,T],\\ \hat s\varphi_*\circ F^{-1}(\frac T2,x_1)\quad\mbox{for}\,\, x_0\in [0,\frac T2],\end{matrix}\right.\end{equation}
 where function $\varphi_*$ defined by (\ref{gopnikz}) and diffeomorphism $F$ is constructed in the beginning of the proof of this proposition, parameter $\hat s>s_0.$

From (\ref{lobster}) and (\ref{main1}) we obtain (\ref{mika1}). Proof of Proposition \label{mika} is complete. $\blacksquare$

{\bf Corollary 1.} {\it Let $f_1\in L^2(Q_+), f_2\in L^2(Q_-), \tilde r\in L^2(0,T)$ and all conditions of Proposition \ref {mika}  holds true. Then there exists function $\eta(x_1)\in C^2(\bar \Omega)$ , $\eta(x_1) <0$ on $\bar\Omega$ and  a constant $C_{80}$ independent of  $v=(v_1,v_2)$ such that
\begin{eqnarray}\label{mika1}
\sum_{\vert\alpha\vert\le 1}\Vert ((T-x_0)^{-3})^\frac {(3-2\vert\alpha\vert)}{2} \partial^\alpha  v_2\,
e^{\psi^*}\Vert_{L^2(Q_-)}  +\sum_{\vert\alpha\vert\le 1}\Vert ((T-x_0)^{-3})^{\frac{5-2\vert\alpha\vert}{2}} \partial^\alpha  v_1\,
e^{\psi^*}\Vert_{L^2(Q_+)}\nonumber\\
+ \Vert ((T-x_0)^{-3})^\frac 32\partial_{x_1}^+v e^{\psi^*}\Vert_{L^2(0,T)}+\Vert ((T-x_0)^{-3})^\frac 12\partial_{x_1}^-v e^{\psi^*}\Vert_{L^2(0,T)}\nonumber\\+\Vert ((T-x_0)^{-3})^\frac 12\partial_{x_0} v e^{\psi^*}\Vert_{L^2(0,T)}+\Vert ((T-x_0)^{-3})^\frac 52 v e^{\psi^*}\Vert_{L^2(0,T)}\nonumber\\
\le C_{80}(\Vert  (T-x_0)^{-3}f_1 e^{\psi^*}\Vert_{L^2(Q_+)}+\Vert  f_2 e^{\psi^*}\Vert_{L^2(Q_-)}\nonumber\\+\Vert (T-x_0)^{-\frac 32}\tilde re^{\psi^*})\Vert_{L^2(0,T)}
+\Vert ((T-x_0)^{-3})^\frac 52v e^{\psi^*}\Vert_{L^2(Q_\omega)}),
\end{eqnarray} where $\psi^*(x)=\eta(x_1)/(T-x_0)^3.$
}

\section{ Proof of Theorem \ref{gonduras}.}

We consider the linearization of the null controllability problem (\ref{eq1}) - (\ref{eq5}):
\begin{equation}\label{Peq1}
L_1(x,D)z_1=\rho_1\partial_{x_0} z_1- a_1\partial_{x_1}^2z_1+b_1\partial_{x_1}z_1+c_1z_1=f_1+\chi_\omega u\quad \mbox{in}\,\, Q_+,
\end{equation}
\begin{equation}
\label{Pe2}
L_2(x,D)z_2=\rho_2\partial_{x_0} z_2-a_2\partial_{x_1}^2z_2+b_2\partial_{x_1}z_2+c_2z_2=f_2\quad \mbox{in}\,\, Q_-,
\end{equation}
\begin{equation}\label{Peq3}
z_1(x_0,0)-z_2(x_0,0)=\partial_{x_1}z_1(x_0,0)-\partial_{x_1}z_2(x_0,0)-M\partial_{x_0}z_1(x_0,0)-r(x_0)=0 \quad \mbox{on} \,\,[0,T],
\end{equation}
\begin{equation}\label{Peeq4}
z_1(x_0,b)=z_2(x_0,a)=0 \quad \mbox{on} \,\,[0,T],
\end{equation}
\begin{equation}\label{Peq4}
z(0,x_1)=z_0(x_1).
\end{equation}
Here
$$
z=\left \{ \begin{matrix} z_1\quad \mbox{for}\quad x\in Q_+,\\ z_2\quad \mbox{for} \quad x\in Q_-.\end{matrix}\right.
$$
  We are looking for control $u$ such that at moment $T$ we have
\begin{equation}\label{Peq5}
z(T,\cdot)=0.
\end{equation}
We have

\begin{proposition}\label{gondurass} Suppose that assumptions (\ref{eq5})- (\ref{eq8}) holds true, $\hat s>s_0$ where $s_0$ is the parameter from Proposition \ref{zoopa}.  Let $(T-x_0)^\frac{15}{2}e^{-\psi^*}f_1\in L^2(Q_+), (T-x_0)^\frac{9}{2}e^{-\psi^*}f_2\in L^2(Q_-), (T-x_0)^\frac{15}{2}e^{-\psi^*}r\in L^2(0,T),$    $z_0\in H^1_0(\Omega)$ and function $\psi^*$ is defined by (\ref{lodka}). Then  controllability problem (\ref{Peq1})-(\ref{Peq5}) has solution $(z,u)\in H^{1,2}(Q_+)\cap H^{1,2}(Q_-)\cap C^0(0,T; H^1(\Omega))\times L^2(Q), \mbox{supp}\, u\subset Q_\omega $ and satisfies the a priori estimate
\begin{eqnarray}\label{vostok}
\sum_{\vert\alpha\vert\le 2}\Vert (T-x_0)^{3\vert\alpha\vert} \partial^\alpha  z\,
e^{-\psi^*}\Vert_{L^2(Q_-)}  +\sum_{\vert\alpha\vert\le 2}\Vert (T-x_0)^{3(\vert\alpha\vert+1)} \partial^\alpha  z\,
e^{-\psi^*}\Vert_{L^2(Q_+)}\nonumber\\
+ \Vert  (T-x_0)^6 z(\cdot,0) e^{-\psi^*}\Vert_{H^1(0,T)}
\le C_1 (\Vert (T-x_0)^\frac{15}{2} fe^{-\psi^*}\Vert_{L^2(Q_+)}\nonumber\\+\Vert (T-x_0)^\frac{9}{2} fe^{-\psi^*}\Vert_{L^2(Q_-)}+\Vert (T-x_0)^\frac{15}{2} r e^{-\psi^*}\Vert_{L^2(0,T)}+\Vert z_0\Vert_{H^1_0(\Omega)}).
\end{eqnarray}
\end{proposition}

{\bf Proof.} Let $\epsilon \in (0,\frac{T}{2})$ be small positive parameter. We set
$\psi^*_\epsilon=\psi^*$ for $x_0\in [0,T-\epsilon ]$ and $\psi^*(x)=\psi^*( T-\epsilon,x_1)$ for $x_0\in [T-\epsilon,T]$ where function $\psi^*$ is constructed in (\ref{lodka}). Denote $\tilde L(x,D)=(\tilde L_1(x,D),\tilde L_2(x,D)), \tilde L_j(x,D) =\frac{1}{a_j} L_j(x,D), $ $\tilde f=f/a_1$ on $Q_+$ and $\tilde f=f/a_2$ on  $Q_-.$
Consider the minimization problem
\begin{eqnarray}\label{A1} J_\epsilon(z,u)=\Vert (T-x_0)^3e^{-\psi^*_\epsilon} z\Vert^2_{L^2(Q_+)}+\Vert e^{-\psi^*_\epsilon} z\Vert^2_{L^2(Q_-)}+ \Vert (T-x_0)^{\frac {15}{2}} e^{-\psi^*}u\Vert^2_{L^2(Q_\omega)}\\+\frac 1\epsilon \Vert  \tilde L(x,D) z-\chi_\omega u-\tilde f\Vert^2_{L^2(Q)}\rightarrow \inf ,\nonumber
\end{eqnarray}
\begin{equation}\label{moloko}
z_1(\cdot,b)=z_2(\cdot,a)=0,\quad  [z](\cdot,0)= [\partial_{x_1} z](\cdot,0)-M\partial_{x_0} z(\cdot,0) +r(\cdot)=0\quad \mbox{on}\,\, [0,T],
\end{equation}
\begin{equation}\label{A2}
z(0,\cdot)=z_0.
\end{equation}
There exists a unique solution to the problem  (\ref{A1})-(\ref{A2}) which we denote as $ (\hat z_\epsilon,\hat u_\epsilon).$ By Theorem \ref{apple} the functions $ (\hat z_\epsilon,\hat u_\epsilon)$ belong to the space $H^{1,2}(Q_+)\cap H^{1,2}(Q_-)\times L^2(Q_\omega).$ Setting $m(x)=1$ for $x\in Q_-$ and $m(x)=(T-x_0)^{15}$ for $x\in Q_+,$ $p_\epsilon = m(x)\frac{e^{-2\psi^*_\epsilon}}{\epsilon}(\tilde L(x,D) \hat z_\epsilon-\chi_\omega\hat u_\epsilon-\tilde f) \in L^2(Q)$ by the Fermat theorem we have
\begin{equation}\label{lincoln}
J_\epsilon'(\hat z_\epsilon,\hat u_\epsilon)[\delta]=0\quad \forall \delta\in\tilde{\mathcal X},
\end{equation} where
$$
\tilde{\mathcal X}=\{\delta=(\delta_1,\delta_2)\vert \,\,\tilde L(x,D) \delta_1=0 \quad\mbox{in}\,\, Q\setminus [0,T]\times \{0\},\,\,\delta_1(0,\cdot)=0, \,\, \delta_1(\cdot,0)=\delta_2(\cdot,0) \quad \mbox{on}\,\, [0,T],
$$
$$ \delta_1(\cdot,b)=\delta_2(\cdot,a)=0,\quad  ([\partial_{x_1} \delta_1]-M\partial_{x_0} \delta_1)(\cdot,0)=0 ,
$$
$$\quad (\delta_1,(T-x_0)^{-\frac {15}{2}} e^{-\psi^*} \delta_2)\in H^{1,2(Q_+)}\cap H^{1,2}(Q_-)\times L^2(Q_\omega)\}
$$ equipped with the norm
$$
\Vert p\Vert_{\tilde{\mathcal X}}=(\Vert(p_1, (T-x_0)^{\frac {15}{2}} e^{- \psi^*} p_2)\Vert^2_{L^2(Q)\times L^2(Q_\omega)}+\Vert p_1\Vert^2_{H^{1,2}(Q_+)}+\Vert p_1\Vert^2_{H^{1,2}(Q_-)})^\frac 12.
$$
From (\ref{lincoln}) we have
\begin{equation}\label{samurai}
\tilde L^*(x,D) p_\epsilon =-e^{-2\psi^*_\epsilon} \hat z_\epsilon\ \quad\mbox{in}\,\, Q\setminus [0,T]\times \{0\},
\end{equation}
\begin{equation}p_{\epsilon,1}(\cdot,b)=p_{\epsilon,2}(\cdot,a)=0,\quad [p_\epsilon](\cdot,0)=([\partial_{x_1} p_\epsilon]-M\partial_{x_0} p_\epsilon)(\cdot,0)=0\quad \mbox{on}\,\, [0,T],
\end{equation}
\begin{equation}\label{samurai1}
p_\epsilon =e^{-2\psi^*}(T-x_0)^{15} \hat u_\epsilon\quad\mbox{on}\,\,Q_\omega,
\quad p_\epsilon(T,\cdot)=0.
\end{equation}
By Theorem \ref{apple} $p_\epsilon\in H^{1,2}(Q_+)\cap H^{1,2}(Q_-).$
By  Corollary 1 the following estimate is  true:
\begin{eqnarray}\label{nona}
\Vert(T-x_0)^\frac{15}{2} e^{\psi^*}p_\epsilon\Vert_{L^2(Q_+)}+\Vert(T-x_0)^\frac{9}{2} e^{\psi^*}p_\epsilon\Vert_{L^2(Q_-)}\\+\Vert(T-x_0)^\frac{15}{2} e^{\psi^*}p_\epsilon(\cdot,0)\Vert_{L^2(0,T)}\nonumber\\\le C_2(\Vert e^{-\psi^*} (T-x_0)^{\frac {15}{2}} \hat u_\epsilon \Vert_{L^2(Q_\omega)}+\Vert (T-x_0)^3 e^{\psi^*-2\psi^*_\epsilon} \hat z_\epsilon\Vert_{L^2(Q_+)}+\Vert e^{\psi^*-2\psi^*_\epsilon} \hat z_\epsilon\Vert_{L^2(Q_-)}).\nonumber
\end{eqnarray}
Taking the scalar product in $L^2(Q)$ of equation (\ref{samurai}) with function $\hat z_\epsilon$ and integrating by parts we have
$$
2J_\epsilon(\hat z_\epsilon,\hat u_\epsilon)= -(p_\epsilon, \tilde f)_{L^2(Q)}-(\frac{p_\epsilon}{a}(0,\cdot), z_0)_{L^2(\Omega)} +(r, p_\epsilon (\cdot,0))_{L^2(0,T)},
$$ where $a(x)=a_1(x)$ on $Q_+$ and $a(x)=a_2(x)$ on $Q_-.$
Hence by (\ref{nona}) we obtain
\begin{eqnarray}\label{milon}
J_\epsilon(\hat z_\epsilon,\hat u_\epsilon)\le C_3 (\Vert (T-x_0)^\frac{15}{2} fe^{-\psi^*}\Vert_{L^2(Q_+)}\\+\Vert (T-x_0)^\frac{9}{2} fe^{-\psi^*}\Vert_{L^2(Q_-)}+\Vert r e^{-\psi^*}(T-x_0)^{\frac{15}{2}}\Vert^2_{L^2(0,T)}+\Vert z_0\Vert^2_{L^2(\Omega)}).\nonumber
\end{eqnarray}
From the sequence $(\hat z_\epsilon, \hat u_\epsilon)$ one can take a subsequence $(e^{-\psi^*_{\epsilon_j}}\hat z_{\epsilon_j}, \hat u_{\epsilon_j})$ which converges weakly to $(e^{-\psi^*}z,u)$ in the space $L^2_m(Q)\times L^2_{(T-x_0)^{15}e^{-2\psi^*}}(Q_\omega).$ From (\ref{milon}) we have
\begin{eqnarray}\label{milon1}
\Vert (T-x_0)^3e^{-\psi^*} z\Vert^2_{L^2(Q_+)}+\Vert e^{-\psi^*} z\Vert^2_{L^2(Q_-)}\le C_4 (\Vert (T-x_0)^\frac{15}{2} fe^{-\psi^*}\Vert_{L^2(Q_+)}\\+\Vert (T-x_0)^\frac{9}{2} fe^{-\psi^*}\Vert_{L^2(Q_-)}+\Vert r e^{-\psi^*}(T-x_0)^{\frac{15}{2}}\Vert_{L^2(0,T)}+\Vert z_0\Vert_{L^2(\Omega)}).\nonumber
\end{eqnarray}
 By Theorem \ref{apple} function $z$ belongs to the space $ H^{1,2}(Q_+)\cap H^{1,2}(Q_-).$
From (\ref{milon1}) and (\ref{lobster})
\begin{eqnarray}\label{milon2}
\sum_{\vert\alpha\vert\le 1}\Vert (T-x_0)^{3\vert\alpha\vert} \partial^\alpha  z\,
e^{-\psi^*}\Vert_{L^2(Q_-)}  +\sum_{\vert\alpha\vert\le 1}\Vert (T-x_0)^{3(\vert\alpha\vert+1)} \partial^\alpha  z\,
e^{-\psi^*}\Vert_{L^2(Q_+)}\nonumber\\
\le C_5 (\Vert (T-x_0)^\frac{15}{2} fe^{-\psi^*}\Vert_{L^2(Q_+)}\nonumber\\+\Vert (T-x_0)^\frac{9}{2} fe^{-\psi^*}\Vert_{L^2(Q_-)}+\Vert  (T-x_0)^\frac{15}{2} r e^{-\psi^*}\Vert_{L^2(0,T)}+\Vert z_0\Vert_{L^2(\Omega)}).
\end{eqnarray}
The the pair $(z,u)$ satisfies  equations (\ref{moloko}), (\ref{A2}). Estimate (\ref{vostok}) follows from (\ref{milon2}) and (\ref{lobster}).
Proof of proposition is complete.
$\blacksquare$

{\bf Proof of Theorem \ref{gonduras}.} We set $\mu_1=(T-x_0)^{15} m(x) e^{-2\psi^*}, m(x)= (T-x_0)^{15}$ for $x\in Q_+$, $m(x)= (T-x_0)^{9}$ for $x\in Q_-,$ $\mu_2=(T-x_0)^{15}e^{-2\psi^*},$ $\mbox{\bf Y}=L^2_{\mu_1}(Q)\times L^2_{\mu_2}(0,T)\times H^1_0(\Omega)$ and $\mbox{\bf X}=\{ w=(w_1,w_2)\vert L_1(x,D)w_1\in L^2_{\mu_1}(Q_+),\,\,  L_1(x,D)w_2\in L^2_{\mu_1}(Q_-),\newline  w_1(\cdot,b)=w_2(\cdot,a)=0, \,\,[w](\cdot,0)=0,\,\, ([\partial_{x_1} w]-M\partial_{x_0} w)(\cdot,0)\in L^2_{\mu_2}(0,T), (T-x_0)^{3(\vert\alpha\vert+1)} \partial^\alpha w\in H^{1,2}(Q_+), (T-x_0)^{3\vert\alpha\vert} \partial^\alpha w\in  H^{1,2}(Q_-)\,\, \forall \vert\alpha\vert\le 2\}.$
Let $(\mbox{\bf w},\mbox{\bf u})$ be the pair of functions defined by Condition 1.
Consider the mapping
$$
F(w,u)=(G(w+\mbox{\bf w})+(\chi_\omega u +\chi_\omega \mbox{\bf u},0), ([\partial_{x_1}w]-M\partial_{x_0} w)(0,\cdot), w(0,\cdot)): \mbox{\bf X}\rightarrow \mbox{\bf Y},
$$ where $G=(G_1,G_2).$
Observe that by Condition 1 $F(0,0)=0.$ Obviously mapping $F\in C^1(\mbox{\bf X},\mbox{\bf Y})$ and $F'(0,0)[\delta_1,\delta_2]=(G'(\mbox{\bf w})[\delta_1]+(\chi_\omega \delta_2,0), ([\partial_{x_1}\delta_1]-M\partial_{x_0} \delta_1)(\cdot,0), \delta_1(0,\cdot)),$  where $G'(\mbox{\bf w})[\delta_1]=(L_1(x,D)\delta_{1,1}+\partial_{\xi_1}g_1(x,\mbox{\bf w}, \partial_{x_1}\mbox{\bf w})\delta_{1,1}+\partial_{\xi_2}g_1(x,\mbox{\bf w}, \partial_{x_1}\mbox{\bf w})\partial_{x_1}\delta_{1,1}, L_2(x,D)\delta_{1,2}+\partial_{\xi_1}g_2(x,\mbox{\bf w}, \partial_{x_1}\mbox{\bf w})\delta_{1,2}+\partial_{\xi_2}g_2(x,\mbox{\bf w}, \partial_{x_1}\mbox{\bf w})\partial_{x_1}\delta_{1,2}). $ Since $\mbox{\bf w}\in L^\infty(Q)$   by (\ref{eq9}) function \newline  $\partial_{\xi_2}g_2(x,\mbox{\bf w}, \partial_{x_1}\mbox{\bf w})\in L^\infty(Q_-)$ and function  $\partial_{\xi_2}g_1(x,\mbox{\bf w}, \partial_{x_1}\mbox{\bf w})\in L^\infty(Q_+)$  and $\partial_{\xi_2}g_1(x,\mbox{\bf w}, \partial_{x_1}\mbox{\bf w})\in L^2(Q_+), \partial_{\xi_1}g_2(x,\mbox{\bf w}, \partial_{x_1}\mbox{\bf w})\in L^2(Q_-).$  Since $\mbox{\bf X}\subset L^\infty(Q)$ operator $G'\in \mathcal L(\mbox{\bf X},\mbox{\bf Y})).$ By Proposition \ref{gondurass} $\mbox{Im}\, F'(0,0)=\mbox{\bf Y}.$ Hence by the implicit function theorem the equation $F(w,u)=\mbox{\bf z}$ can be solved for all $\mbox{\bf z}\in \mbox{\bf X}$ from some neighborhood of $0$ in the space $\mbox{\bf Y}.$ Proof of the theorem is complete.
$\blacksquare$

\section{ APPENDIX}

\begin{theorem}\label{apple} Let $u\equiv 0, g_1\equiv 0,g_2\equiv 0.$  Suppose that (\ref{eq6})-(\ref{eq7}) holds true.
Then for any $w_0\in H^1_0(\Omega),$ $f_1\in L^2(Q_+)$, $f_2\in L^2(Q_-), \tilde r\in L^2(0,T)$ there exist a unique solution to problem (\ref{eq1})-(\ref{eq4}) such that
\begin{eqnarray}\label{lobster}
\Vert w_1\Vert_{H^{1,2}(Q_+)}+\Vert w_2\Vert_{H^{1,2}(Q_-)}+\Vert w_1(\cdot,0)\Vert_{H^1(0,T)}
\\\le C_1(\Vert f_1\Vert_{L^2(Q_+)}+\Vert f_2\Vert_{L^2(Q_-)}+\Vert w_0\Vert_{H^1(\Omega)}+\Vert \tilde r\Vert_{L^2(0,T)}).\nonumber
\end{eqnarray}
\end{theorem}

{\bf Proof.}  Without loss of generality we may assume that $a_j\equiv 1.$ First we observe that it suffices to prove Theorem \ref{apple} for the case $w_0\equiv 0$. Indeed consider the function $v_* =w_0-\eta w_0(0)$, where $\eta\in C^\infty_0(\Omega), \eta(0)=1.$ Let $v_1\in H^{1,2}(Q_+)$ and $v_2\in H^{1,2}(Q_-)$ be solution to the following initial value problems
$$
L_1(x,D)v_1=f_1\quad\mbox{in}\,\, Q_+,\quad  v_1(x_0,0)=v_1(x_0,b)=0, \quad v_1(0,\cdot)=v_*,
$$
$$
L_2(x,D)v_2=f_2\quad\mbox{in}\,\, Q_-,\quad  v_2(x_0,a)=v_2(x_0,0)=0, \quad v_2(0,\cdot)=v_*.
$$
Then if $z_1,z_2$ be solution to the problem
\begin{eqnarray}\label{Sskazka1}
L_1(x,D)z_1=g_1\quad \mbox{in}\, \,Q_+,\quad  L_2(x,D)z_2=g_2 \quad \mbox{in}\, \,Q_-,\\
 z_1(0,\cdot)=0,\quad z_2(0,\cdot)=0, \quad z_2(\cdot,a)=z_1(\cdot,b)=0,\nonumber\\ z_1(\cdot,0)=z_2(\cdot,0), \quad \partial_{x_1}z_1-\partial_{x_1}z_2=M\partial_{x_0}z(\cdot,0)+\tilde r\quad \mbox{on}\,\, [0,T],\nonumber
\end{eqnarray}
where $g_1=-L_1(x,D)(\eta w_{0}(0)), g_2=-L_2(x,D)(\eta w_{0}(0)).$ Then the function $(w_1,w_2)=(v_1,v_2)-(\eta w_0(0),\eta w_{0}(0))$ be solution to problem (\ref{eq1})-(\ref{eq4}).
Let $K$ be a large positive parameter. In problem (\ref{eq1})-(\ref{eq4}) we move from unknown function $(w_1,w_2)$ to the unknown function $z=(z_1,z_2)=(w_1,w_2)e^{-Kx_0}.$ the function $z$ satisfies
\begin{eqnarray} \label{Sskazka1}
(L_1(x,D)+K)z_1=g_1\quad \mbox{in}\, \,Q_+,\quad ( L_2(x,D)+K)z_2=g_2 \quad \mbox{in}\, \,Q_-,\\
 z_1(0,\cdot)=0,\quad z_2(0,\cdot)=0,\quad z_2(\cdot,a)=z_1(\cdot,b)=0,\nonumber\\ z_1(\cdot,0)-z_2(\cdot,0)= (\partial_{x_1}z_1-\partial_{x_1}z_2-M\partial_{x_0}z-MKz)(\cdot,0)+p(\cdot)=0\quad \mbox{on}\,\, [0,T],\nonumber
\end{eqnarray}
where
$g= (g_1,g_2), g_i=e^{-Kx_0}f_i, p=e^{-Kx_0}\tilde r.$

Assume that $z\in H^{1,2}(Q_+)\cap  H^{1,2}(Q_-) $ be solution to problem (\ref{Sskazka1}). The there exists a constant $C_2$ independent of $z$ such that
\begin{equation}\label{ispolin}
\Vert z \Vert_{H^{1,2}(Q_+)\cap  H^{1,2}(Q_-)}+\Vert z(\cdot,0)\Vert_{H^1(0,T)}\le C_2(\Vert p\Vert_{L^2(0,T)}+\Vert g\Vert_{L^2(Q)}).
\end{equation}
The estimate (\ref{ispolin}) is the standard energy estimate which can be proved by multiplying  (\ref{Sskazka1}) by $z$ and $\partial_{x_0}z.$
in particular the estimate (\ref{ispolin}) implies that it suffices to prove the existence of solution to problem (\ref{Sskazka1})  assuming that $g\in C^0([0,T]; L^2(\Omega)), p\in C^0[0,T]$ and coefficients $b,c\in C^1(\bar Q).$

Now consider the discretization of problem (\ref{eq1})-(\ref{eq4}) in time
\begin{eqnarray}\label{robespier} \rho (x_{0,k+1},\cdot)\frac{(z_{k+1}-z_{k})}{h}-\partial_{x_1}^2z_{k+1}+ b(x_{0,k+1},\cdot)\partial_{x_1}z_{k+1}\nonumber \\+(c(x_{0,k+1},\cdot)+K)z_{k+1}=g(x_{0,k},\cdot)\quad\mbox{in}\quad \Omega,
\end{eqnarray}
\begin{equation}\label{robespier1}
[z_{k}]=0,\quad\partial_{x_1}[z_{k+1}]=M\frac{(z_{k+1}-z_{k})(0)}{h}+ MKz_{k+1}+p(x_{0,k}).
\end{equation}
Here $z_k=(z_{k,1},z_{k,2}), z_0(x_1)=0$, $x_{0,k}=kT/N, h=\frac{T}{N}, N\in \Bbb Z_+.$

Multiplying equation (\ref{robespier}) by $h z_{k+1},$ integrating by parts and taking sum respect to $k$ we have
\begin{eqnarray}\label{grom0}
\sum_{k=0}^{N-1}\left\{ \int_\Omega\rho (x_{0,k+1},x_1)(z_{k+1}-z_{k})z_{k+1}dx_1+ M\left(\frac{z_{k+1}-z_{k}}{h}\right)(0)z_{k+1}(0)\right. \nonumber\\ \left.+ KMz_{k+1}^2(0) +h(\Vert\partial_{x_1} z_{k+1}\Vert^2_{L^2(\Omega)}+K\Vert z_{k+1}\Vert^2_{L^2(\Omega)})\right\}=
\\ \nonumber
\sum_{k=0}^{N-1}(g(x_{0,k},\cdot)- b(x_{0,k+1},\cdot)\partial_{x_1}z_{k+1}- c(x_{0,k+1},\cdot)z_{k+1},z_{k+1})_{L^2(\Omega)}-p(x_{0,k})z_{k+1}(0)\le
\\\nonumber
C_3\sum_{k=0}^{N-1}\left\{\Vert g(x_{0,k},\cdot)\Vert^2_{L^2(\Omega)}+\frac{M}{2}\Vert\partial_{x_1}z_{k+1}\Vert^2_{L^2(\Omega)}+ p^2(x_{0,k})\right\} +\frac{K}{2}\sum_{k=0}^{N-1}\Vert z_{k+1}\Vert^2_{L^2(\Omega)}.
\end{eqnarray}
Observe that
\begin{eqnarray}\label{grom1}
\sum_{k=0}^{N-1}\int_\Omega\rho (x_{0,k+1},x_1)(z_{k+1}-z_{k})z_{k+1}dx_1=
\\ \nonumber \sum_{k=0}^{N-1}\int_\Omega\left (\rho (x_{0,k+1},x_1)z^2_{k+1}-\root\of{\rho (x_{0,k},x_1)}z_{k}\root\of{\rho (x_{0,k+1},x_1)}z_{k+1}\right)dx_1
\\ \nonumber
+\int_{\Omega}\left (\frac{\root\of{\rho (x_{0,k},x_1)}-\root\of{\rho (x_{0,k+1},x_1)}}{h}\right ) hz_kz_{k+1}\,\root\of{\rho (x_{0,k+1},x_1)}dx_1\ge
\\ \nonumber
-C_4\Vert \partial_{x_0}\tilde\rho\Vert_{L^\infty(Q)}\sum_{k=0}^{N} \Vert z_k\Vert^2_{L^2(\Omega)}.
\end{eqnarray}
From (\ref{grom0}), (\ref{grom1}) we obtain
\begin{equation} \label{zampolit}
\sum_{k=0}^{N} (K\Vert z_k\Vert^2_{L^2(\Omega)}+ \Vert\partial_{x_1} z_k\Vert^2_{L^2(\Omega)})\le C_5(\Vert g\Vert_{C^0([0,T];L^2(\Omega))} +\Vert p\Vert_{C^0[0,T]}).
\end{equation}
Multiplying equation (\ref{robespier}) by $\frac{(z_{k+1}-z_{k})}{h},$ integrating by parts and taking sum respect to $k$ we have
$$
\sum_{k=0}^{N_0-1}\left\{ \int_\Omega\rho (x_{0,k},x_1)\left (\frac{z_{k+1}-z_{k}}{h}\right )^2dx_1+ M\left(\frac{z_{k+1}-z_{k}}{h}\right)^2(0)+\frac 1 h\Vert\partial_{x_1} z_{k+1}\Vert^2_{L^2(\Omega)}\right.
$$
$$\left.+(MKz_{k+1},\frac{(z_{k+1}-z_{k})}{h})_{L^2(\Omega)}
-\frac 1h(\partial_{x_1}z_{k+1},\partial_{x_1}z_k)_{L^2(\Omega)}\right\}=
$$
$$
\sum_{k=0}^{N_0-1}(g(x_{0,k},\cdot)),\frac{(z_{k+1}-z_{k})}{h})_{L^2(\Omega)}+ p(x_{0,k})\frac{(z_{k+1}-z_{k})}{h}\le
$$
$$
C_6\sum_{k=0}^{N_0-1}\left\{\Vert g(x_{0,k},\cdot)\Vert^2_{L^2(\Omega)}+\Vert\partial_{x_1}z_{k+1}\Vert^2_{L^2(\Omega)}+p^2(x_{0,k})+\Vert z_{k+1}\Vert^2_{L^2(\Omega)}\right\}
$$
$$+\sum_{k=0}^{N_0-1}\left (\frac \alpha 2\Vert \frac{z_{k+1}-z_{k}}{h}\Vert^2_{L^2(\Omega)}+\frac M2 \left (\frac{(z_{k+1}-z_{k})}{h}\right )^2\right ).
$$
This inequality and (\ref{zampolit}) imply
\begin{eqnarray}\label{rom}
\sum_{k=0}^{N-1}\left\{ h\int_\Omega\left (\frac{z_{k+1}-z_{k}}{h}\right )^2dx_1+ M h\left(\frac{z_{k+1}-z_{k}}{h}\right)^2(0)+\frac 1 h
\mbox{sup}_{k\in \{1,\dots, N\}}\Vert\partial_{x_1} z_{k}\Vert^2_{L^2(\Omega)}\right\}\nonumber\\
\le C_7\sum_{k=0}^{N-1}h\{ \Vert g(x_{0,k},\cdot)\Vert^2_{L^2(\Omega)}+p^2(x_{0,k})\}.
\end{eqnarray}
We define function $\tilde z_N$ as follows: $\tilde z_N(x_0,\cdot)=z_k$  if $x_0=x_{0,k}$, otherwise on interval $(T(k+1)/N,Tk/N)$ function $\tilde z_N$ is the linear function.
Then estimate (\ref{rom}) implies
\begin{equation}\label{grom}
\Vert \tilde z_N\Vert_{H^{1,2}(Q)}+  K^\frac 14\Vert \tilde z_N\Vert_{L^2(0,T;H^1(\Omega))}+\root\of {K}\Vert\tilde z_N\Vert_{L^2(Q)}\le C_8(\Vert p\Vert_{C^0[0,T]}+\Vert g\Vert_{C^0([0,T];L^2(\Omega))}).
\end{equation}
Functions $\tilde z_N=(\tilde z_{N,1}, \tilde z_{N,2})$ satisfy the initial value problem
\begin{eqnarray} \label{skazka11}
(L_1(x,D)+K)\tilde z_{N,1}=g_{1,N}\quad \mbox{in}\, \,Q_+,\quad ( L_2(x,D)+K)\tilde z_{N,2}=g_{2,N} \quad \mbox{in}\, \,Q_-,\\
 \tilde z_{N,1}(0,\cdot)=0,\quad \tilde z_{N,2}(0,\cdot)=0,\quad \tilde z_{N,2}(\cdot,a)=\tilde z_{N,1}(\cdot,b)=0,\nonumber\\ \tilde z_{N,1}(\cdot,0)-\tilde z_{N,2}(\cdot,0)=(\partial_{x_1}\tilde z_{N,1}-\partial_{x_1}\tilde z_{N,2}-M\partial_{x_0}\tilde z_{N,1}(\cdot,0)-KM\tilde z_{N,1})(\cdot,0)+p_N (\cdot)=0\,\,\mbox{on}\,\, [0,T].\nonumber
\end{eqnarray}

Let  sequence $\tilde z_{N}$ after possibly taking a subsequence converges to the function $z$ weakly in  $H^{1,2}(Q_+)\cap H^{1,2}(Q_-).$
Observe that
$$
(g_{1,N},g_{2,N}, p_N)\rightarrow 0\quad \mbox{weakly} \,\, \mbox{in}\, L^2(0,T;H^{-1}(\Omega_+))\times L^2(0,T;H^{-1}(\Omega_-))\times H^{-1}(0,T).
$$
Passing to the limit in (\ref{skazka11}) we obtain that function $z$ solution to problem (\ref{Sskazka1}). Estimate (\ref{lobster}) follows from (\ref{ispolin}). Proof of theorem is  complete.
$\blacksquare$


\begin{thebibliography}{99} %

\bibitem{AB} J. Ben Amara, H. Bouzidi, {\it Null boundary controllability of a one-dimensional heat equation with an internal point mass and variable coefficients,}
J. Math. Physics, {\bf 59}, (2018), 011512.

\bibitem{AE2018}    S. Avdonin, J. Edwards, {\it Exact controllability for string with attached masses,}
SIAM J.\ Control Optim., {\bf  56}, (2018), 011512.

\bibitem{AYR} A. Benabdallah, Y. Dermenjiian and  J. L. Rousseau {\it Carleman estimates for the one-dimensional heat equation with a discontinuous coefficient and applications  to controllability and an inverse problem,} J. Math. Analysis and Applications, {\bf 336}, (2007), 865-887.

    \bibitem{Dubova} A. Doubova, A. Osses and J.-P. Puel, {\it  Exact controllability to trajectories for the semilinear heat equation with discontinuous coefficients,}
ESAIM COCV., {\bf 8}, (2002),  621-661.
%
\bibitem{D-F} A. Dubova, E. Fernander-Cara, {\it Some control results for simplified one-dimensional models of fluid-solid interactions,} Mathematical Models and Methods in Applied Sciences, {\bf 15}, (2005), 783-824.

\bibitem{HM1} S. Hansen, J. Martinez,
{\it Modeling of a heat equation with Dirac density, }
Proceedings of Dynamic Systems and Applications, {\bf 7}, (2016),  148-154.

\bibitem{HM2} S. Hansen,  J. Martinez,
{\it Null boundary controllability of a 1-dimensional heat equation with an internal point mass,}
Decision and Control (CDC), 2016 IEEE 55th Conference,  (2016), 4803-4808.


\bibitem{Scott} S. Hansen,  E Zuazua,
{\it Exact controllability and stabilization of a vibrating string with an interior point mass,}
SIAM J. Cont. Optim., {\bf 33}, (1995), 1357-1391.

\bibitem{H}
L. H\"ormander, {\it Linear Partial Differential Operators},
Spring-Verlag, Berlin, 1963.


\bibitem{Im}
O. Imanuvilov, {\it
Controllability of parabolic equations,} Math. Sb., {\bf 186}, (1995),
879-900.

\bibitem{IM}
O. Imanuvilov, M. Yamamoto, {\it
Carleman estimate and an inverse source problem for the Kelvin-Voigt model for
viscoelasticity,} arXiv:1711.09276 .

\bibitem{LTT} Y. Liu, T. Takahashi and M. Tucsnak, {\it Single input controllability  of a simplified fluid-structure model,} ESIAM: COCV., {\bf 19}, (2013), 20-42.
\bibitem{MTT} D. Maity, T. Takahashi  and M. Tucsnak, {\it Analysis of a system modeling of a piston in a viscous gas,} J. Math. Fluid Mechanics, {\bf 19}, (2017), 551-579.

    \bibitem{Russeau} J. L. Russeau, M. Leautaud and  L. Robbiano, {\it Controllability of parabolic system with diffusive interface,} J. Eur. Math. Soc., {\bf 15}, (2013), 1485-1574.

\bibitem{Sog} C. D. Sogge,
{\it Fourier Integrals in Classical Analysis,}
Cambridge University Press, Cambridge, 1993.

\bibitem{Stein} E. M. Stein, {\it Harmonic Analysis: Real-Variable Methods,
Orthogonality and Oscillatory Integrals,} Princeton University Press,
Princeton, New Jersey, 1993.

\bibitem{Tay} M. Taylor, {\it Pseudodifferential Operators and Nonlinear PDE,}
Birkh\"auser, Basel, 1991.

\bibitem{Tay1} M. Taylor, {\it Tools for PDE: Pseudodifferential Operators,
Paradifferential Operators, and Layer Potentials,}
American Mathematical Society, Rhode Island, 2000.

\bibitem{VZ} J.L. Vazquez, E. Zuazua, {\it Large behavior for a simplified  1D model of fluid-solid iteraction,} Comm. P.D.E.,  {\bf 28}, (2003), 1705-1738.

\end{thebibliography}
\end{document}